 \documentclass[bj]{imsart}

 \RequirePackage[OT1]{fontenc}
 \RequirePackage{amsthm,amsmath,natbib}
 \RequirePackage[colorlinks,citecolor=blue,urlcolor=blue]{hyperref}

 \usepackage[top=4.38cm,bottom=4.4cm, outer=3.5cm, inner=3.5cm]{geometry}

 \usepackage{tikz} 
 \usepackage{amsmath}
 \usepackage{mathrsfs} 
 \usepackage{amsfonts} 
 \usepackage{pgfplots}

 \usepackage{amssymb}

 \def\ar{\!\!\!&} 

 \arxiv{arXiv:1912.04795}

 \startlocaldefs
 \numberwithin{equation}{section}
 \theoremstyle{plain}
 \newtheorem{theorem}{Theorem}[section]
 \newtheorem{condition}[theorem]{Condition}
 
 \newtheorem{lemma}[theorem]{Lemma}
 
 \newtheorem{proposition}[theorem]{Proposition}
 \theoremstyle{remark}
 \newtheorem{remark}[theorem]{Remark}

 \endlocaldefs

\begin{document}

 \begin{frontmatter}
 \title{Asymptotic Results for Heavy-tailed L\'evy Processes and their Exponential Functionals}
 \runtitle{L\'evy Processes and their Exponential Functionals}

 \begin{aug}
 \author{\fnms{Wei} \snm{Xu}\thanksref{a,e1}
	\ead[label=e1,mark]{xuwei@math.hu-berlin.de}}

 \address[a]{Department of Mathematics, Humboldt-Universit\"at zu Berlin, Unter den Linden 6, 10099 Berlin, Germany.
 \printead{e1}}

\runauthor{W. Xu}

\affiliation{Humboldt-Universit\"at zu Berlin}

\end{aug}

 \begin{abstract}
 In this paper we first provide several conditional limit theorems for L\'evy processes with negative drift and regularly varying tail.
 Then we apply them to study the asymptotic behavior of expectations of some exponential functionals of heavy-tailed L\'evy processes.
 As the key point, we observe that the asymptotic mainly depends on the sample paths with early arrival of large jump.
 Both the polynomial decay rate and the exact expression of the limit coefficients are given.
 As an application, we give an exact description for the extinction speed of continuous-state branching processes in heavy-tailed L\'evy random environment with stable branching mechanism.
 \end{abstract}

 \end{frontmatter}

 \section{Introduction}
 \setcounter{equation}{0}
 \medskip

 The long-run behavior of L\'evy processes and their functionals has been extensively explored in the literature in the past decades.
 This is mainly justified by their wide and important applications in various fields such as risk theory, mathematical finance, physics and population evolution.
 In this work, we mainly explore the asymptotic behavior of L\'evy processes and their exponential functionals in the presence of regularly varying tails.

 In the first part of this work, we study the asymptotic behavior of heavy-tailed L\'evy process  $\{\xi_t:t\geq 0  \}$ with negative drift conditioned to stay positive, which will play an important role to study their exponential functionals.
 Denote by $\tau_0$ the first entrance time of $\xi$ in $(-\infty,0)$.
 For any $x,t>0$, let $\mathcal{N}_t^x$ be the number of jumps of $\xi$ larger than $x$ up to time $t$ and $\mathcal{J}^x := \inf\{t>0: \mathcal{N}_t^x>0 \}$ be the arrival time of the first one.
 We first show that $\xi_t>0$ for large $t$ if and only if a large jump occurs at a random time that is uniformly distribution in the time  interval $[0,t]$, i.e.
  \begin{eqnarray*}
  \mathbf{P}\big\{\mathcal{N}_t^{at}=1\,\big|\,\xi_t> 0\big\}\sim \mathbf{P}\big\{\xi_t> 0\,\big|\,\mathcal{N}_t^{at}=1\big\}\to 1
  \quad \mbox{and}\quad
  \big\{\mathcal{J}^{at}/t\,\big|\,\xi_t> 0\big\} \to\mathbf{U},
 \end{eqnarray*}
 where $\mathbf{U}$ is a uniformly distributed random variable  on $[0,1]$.

 By contrast, we also show that, in some sense, the process $\xi$ stays positive for long time ($\tau_0>t$ for large $t$) if and only if  a large jump occurs before time $t$ as well as the process needs to stay positive before the arrival of large jump, i.e.
 \begin{equation*}
\mathbf{P}_x\big\{\tau_0> \mathcal{J}^{at}\,\big|\,\tau_0> t\big\}
\sim \mathbf{P}_x\big\{\tau_0> t\,\big|\,\tau_0> \mathcal{J}^{at}\big\} \to 1.
\end{equation*}
  This motivates us to identify that spatial-scaled process $\{ \xi_s/t,\, s\geq 0\,|\, \tau_0>t\}$ converges weakly to a non-degenerate limit, whose sample paths are step functions with a single jump occurring at random time which has a size-biased distribution, i.e.
  \begin{equation*}
  \{t^{-1}\xi_s,s\geq 0|\tau_0>t, \xi_0=x\}\to\{\mathcal{P}_{a,\alpha}\cdot \mathbf{1}_{\{\mathcal {T}_x<s\}},s\geq 0\}
  \end{equation*}
 weakly, where $\mathcal{P}_{a,\alpha}$ and $\mathcal {T}_x $ are two independent positive  random variables.
  Some analogous functional limit theorems for random walk $\{Z_n:n=0,1,\cdots\}$  conditioned to stay positive have been studied by many authors. Here we provide a brief summary.
 Let $\tau_0^d$ denote its first passage time in $(-\infty,0)$.
 For random walk with negative drift and regularly varying tail,  \cite{Durrett1980} proved, in the special case when $Z_1$ has finite variance, that both of the time-spatial rescaled processes
 $\{ Z_{[nt]}/n,\, t\in[0,1]\,|\,Z_n>0\}$ and $\{ Z_{[nt]}/n,\, t\in[0,1]\,|\,\tau^d_0>n\}$
 converge weakly to a non-degenerate limit with sample paths having a single jump at time $0$ and decreasing linearly. These have been partially extended by \cite{DoneyJones2012} to L\'evy processes with infinite variance.
  It is clear that the conditional limit theorems in \cite{Durrett1980} are different to ours.
 For the oscillating case, under the Spitzer's condition \cite{AGKV2005,Doney1985} and \cite{Durrett1978} showed that for some regularly varying function $b_n$ with index $\alpha\in(0,2]$  the rescaled process $\{b^{-1}_n Z_{[nt]},\, t\in[0,1] \,|\, \tau_0^d>n \}$ converges weakly to the meander of a strictly $\alpha$-stable process.

 At the end of the first part of this paper, based on the Wiener-Hopf factorization we also provide exact expressions for the decay rates of the Laplace transforms $\mathbf{E}[e^{-\lambda \xi_t},\xi_t>0]$ and $\mathbf{E}_x[e^{-\lambda \xi_t},\tau_0>t]$, which are applied in the following analysis of the asymptotic behaviors of exponential functionals of L\'evy processes.  In detail, we have
  \begin{equation*}
 \mathbf{E}[e^{-\lambda\xi_t};\xi_t\geq 0]\sim  \frac{\alpha}{a\lambda} \cdot\mathbf{P}\{\xi_1> at\}\quad\mbox{and}\quad
\mathbf{E}_x[e^{-\lambda \xi_t},\tau_0>t] \sim C(x)\cdot \frac{\mathbf{P}\{\xi_1> at\}}{at}.
 \end{equation*}
 The analogous results for random walk have been established in \cite{DenisovVatutinWachtel2014}.

 With the help of  conditional limit theorems of L\'evy processes,
  as a second main contribution we analyze the long-term behavior of their exponential functionals
 \begin{equation}\label{eq1.1}
 A_t(\xi) := \int_0^t e^{- \xi_s}ds, \qquad 0\le t\leq \infty.
 \end{equation}
 Here $A_\infty(\xi)= \infty$ is allowed.
 The study of exponential functionals  has drawn the attention of many researchers because of the considerable role they play in mathematical finance, statistics physics and population evolution in random environment; see \cite{GemanYor1993,ComtetMonthusYor1998,BansayePardoSmadi2013}
 and \cite{CarmonaPetitYor1997}.
 The literature can be roughly divided into two classes: finite-time horizon ($t<\infty$)  and infinite-time horizon ($t=\infty$).

 On the finite-time horizon, via the analysis of the Bernstein-gamma function which represents the Mellin transform of $A_t(\xi)$, \cite{PatieSavov2018} provided necessary and sufficient conditions for the finiteness of its negative and positive moments.
 Recently, \cite{BarkerSavov2019} introduced the bivariate Bernstein-gamma function and applied them to study $A_t(\xi)$ with $\xi$ being a subordinator. They provided an explicit infinite convolution formula for the Mellin transform of $A_t(\xi)$.

 On the infinite-time horizon, a necessary and sufficient condition for $A_\infty(\xi)<\infty$ a.s.  was given in \cite{BertoinYor2005}, i.e., $A_\infty(\xi)<\infty$ if and only if $\xi$ drifts to infinity.
 In this case,  the characterizations and Wiener-Hopf type factorization of the law of $A_\infty(\xi)$ can be found in \cite{CarmonaPetitYor1997,PardoPatieSavov2012,PatieSavov2012} and \cite{Yor1992}.
 For more interesting results and properties of $A_\infty(\xi)$, readers may refer to \cite{BarkerSavov2019,PatieSavov2018} and \cite{Vechambre2019}.
 For the general case with $\xi$ replaced by a non-decreasing Markov additive process, \cite{Stephenson2018} provided several results for $A_\infty(\xi)$ including the explicit formulas for its positive moments and a necessary condition for the finiteness of its exponential moments.
 In the case of $A_\infty(\xi) =\infty$, we are usually interested in the asymptotic behavior of  $F(A_t(\xi))$ for some positive,  decreasing function $F$ on $(0,\infty)$ that vanishes at $\infty$.
 Especially,  much attention has been drawn to the decay rate of the following expectation as $t\to \infty$:
 \begin{equation}\label{eq1.2}
 \mathbf{E}[F(A_t(\xi))]
 =
 \mathbf{E}\Big[F\Big(\int_0^t e^{- \xi_s}ds\Big)\Big],
 \end{equation}
 which is closely related to the long-term properties of random processes in random environment; see \cite{CarmonaPetitYor1997} and \cite{KawazuTanaka1993}.

 To the best of our knowledge, almost all previous works, except \cite{PatieSavov2018}, usually considered the expectation (\ref{eq1.2}) with $\xi $ satisfying some exponential moment condition, e.g., the Laplace exponent $\phi(\lambda):=\log\mathbf{E}[\exp\{\lambda\xi_1\}]<\infty$ for some $\lambda>0$.
 Here we give a brief summary for them.
 Readers may refer to \cite{LiXu2018,PalauPardoSmadi2016} and references therein for details.
 For L\'evy processes with bounded variation,
 \cite{CarmonaPetitYor1997} provided a precise asymptotic behavior for  $\mathbf{E}[(A_t(\xi))^{-p}]$ for some constant $p>0$ satisfying that $\phi'(p)>0$.
 Applying the discretization technique and the asymptotic results proved in \cite{GuivarchLiu2001}, \cite{BansayePardoSmadi2013} provided four regimes for the decay rate of the expectation (\ref{eq1.2}) with $\xi$ being a compound Poisson process and $F(z)\sim Cz^{-p}$. Their approach was extended by \cite{PalauPardoSmadi2016}.
 By observing that the asymptotic only depends on the sample paths with slowly decreasing local infimum, \cite{LiXu2018} provided not only four different regimes for the convergence rate of the expectation (\ref{eq1.2}) but also the exact expressions of the limiting coefficients.
 Their proofs heavily rely on the fluctuation theory and limit theorems of L\'evy processes conditioned to stay positive.

  In the second part of this work, we study the asymptotic behavior of the expectation (\ref{eq1.2}) in the case where $\xi$ has negative drift and regularly varying tail.
  In this case, the exponential moment condition can not be satisfied, i.e., $\mathbf{E}[e^{\lambda\xi_1}]=\infty$ for any $\lambda>0$.
 To our knowledge, \cite{PatieSavov2018} is the only one that provided an exact description for the asymptotic behavior of the expectation (\ref{eq1.2}) with $\xi$ oscillating and satisfying the Spitzer's condition.
  Inspired by the analysis in \cite{LiXu2018} and \cite{VatutinZheng2012},
  we observe that the main contribution to the expectation (\ref{eq1.2}) is given by the scenario of  slowly decreasing local infimum and these sample paths can be identified according to the arrival time of their first large jumps.
 This  helps us to establish the polynomial decay for the expectation (\ref{eq1.2}),  i.e.
 \begin{equation*}
 \mathbf{E}[F(A_t(\xi))] \sim C_F\cdot \mathbf{P}\{\xi_1>at\}.
 \end{equation*}
 Here we also provide the exact expression for the  coefficient $C_F$ by  the conditional limit theorems of the functional of L\'evy process starting from a large jump.

 To illustrate the strength of our asymptotic results for the expectation (\ref{eq1.2}), in the last part of this work we study the long-run behavior of  continuous-state branching processes in heavy-tailed  random environment.
 Let $\{Z^\gamma_t: t\ge 0\}$ be a spectrally positive $(\gamma+1)$-stable process with $\gamma\in(0,1]$.
 Note that $Z^\gamma$ is a Brownian motion when $\gamma=1$.
 For any $x\ge 0$, we consider the unique strong solution to the following stochastic integral equation:
 \begin{equation}\label{eq1.3}
 X_t(x) = x + \int_0^t \sqrt[\gamma+1]{(\gamma+1)c X_{s-}(x) } dZ^\gamma_s+ \int_0^t  X_{s-}(x) dZ^e_s.
 \end{equation}
 where $c\ge 0$ and $\{Z^e_t :  t\ge 0\}$ a L\'{e}vy process with no jump less than $-1$ and independent of $Z^\gamma$.
 The solution is called a \textit{continuous-state branching process in random environment (CBRE-process)} with \textit{stable branching mechanism} and \textit{L\'evy random environment} $Z^e$; see \cite{HeLiXu2018,PalauPardo2018}.
 The construction of CBRE-processes as the scaling limit of rescaled  Galton-Watson processes in random environment (GWRE-processes) can be found in \cite{BansayeCaballeroMeleard2019} and \cite{Kurtz1978}.
 From Section~4 in \cite{HeLiXu2018},  there is another L\'{e}vy process $\{\xi_t: t\ge 0\}$ determined by the environment so that the \textit{survival probability} of the CBRE-process at time $t\ge 0$ is given by
 \begin{equation}\label{eq1.4}
 \mathbf{P}\{X_t(x)>0\} = \mathbf{E}\Big[1-\exp\big\{-x\big(c\gamma  A_t(\gamma\xi)\big)^{-1/\gamma}\big\}\Big],
 \end{equation}
 which clearly is a special case of (\ref{eq1.2}).
 Compared to the abundant literature about survival probabilities of GWRE-processes; e.g. \cite{AGKV2005,BansayeVatutin2017,GeigerKerstingVatutin2003,
 GuivarchLiu2001, VatutinZheng2012} and references therein,
 there can be found only several works about survival probabilities of CBRE-processes.
 With the expression (\ref{eq1.4}),
 the asymptotic of survival probability was studied in  \cite{BansayePardoSmadi2013,BoinghoffHutzenthaler2012,LiXu2018, PalauPardo2015} and \cite{PalauPardoSmadi2016} with $\xi$ satisfying some exponential moment condition.
 Recently,  \cite{BansayePardoSmadi2019} studied the extinction speed of CBRE-processes with general branching mechanism and oscillating environment $\xi$. They provided the exact expression for the extinction rate under a necessary condition that is $\mathbf{E}[
 e^{\theta^+\xi_1}]<\infty$ for some $\theta^+>1$. This  excludes the CBRE-processes with the right-tail of random environment being regularly varying.
 Here we provide an exact expression for the asymptotic behavior of survival probabilities of CBRE-processes with stable branching mechanism and  heavy-tailed environment.
 In detail, the survival probabilities decrease to $0$ at a polynomial rate.
 This is in sharp contrast with the exponential decay of  survival probabilities of CB-processes in light-tailed L\'evy random environment; see \cite{LiXu2018}.

 The remainder of this paper is organized as follows.
 In Section~2, we recall definitions and some basic elements of fluctuation theory for L\'evy processes, which are necessary for the proofs.
 In Section~3, we   give several conditional limit theorems for L\'evy processes  with negative drift and regularly varying tail.
 The exact expression of asymptotic behavior of the expectation (\ref{eq1.2}) is given in Section~4.
 In Section~5, we study the decay rate of survival probabilities of CBRE-processes.

 \medskip

 \noindent \textbf{Notation.}
 For any two sets $A,B$, let $A \Delta B $ be their symmetric difference, i.e., $A \Delta B:= (A \cap B^{\rm c}) \cup (B\cap A^{\rm c})$.
 For any $x\in\mathbb{R}$,  let $[x]$ be its integer part and $x^+ := x\vee 0$.
 We make the convention that for any $t_1\leq t_2\in\mathbb{R}$,
 \begin{equation*}
 \int_{t_1}^{t_2}=\int_{(t_1,t_2]}\quad \mbox{and} \quad \int_{t_1}^{\infty}=\int_{(t_1,\infty)}.
 \end{equation*}


 \section{Preliminaries}
 \setcounter{equation}{0}

 In this section, we recall some basic notation and elements of fluctuation theory for L\'evy processes.
 The reader may refer to \cite{Bertoin1996} and \cite{Kyprianou2006} for details.
 Let $(\Omega,\mathscr{F},\mathbf{P})$ be a complete probability space endowed with filtration $\{\mathscr{F}_t\}_{t\geq 0}$ satisfying the usual hypotheses and $\mathbf{D}([0,\infty),\mathbb{R})$ the space of c\'adl\'ag paths endowed with Skorokhod topology.
 Denote by $\xi = \{\xi_t: t\ge 0\}$ a one-dimensional L\'evy process with characteristic exponent ${\it\Phi}(\lambda):=-\log\mathbf{E}[\exp\{\mathtt{i}\lambda\xi_1\}]$. Here ${\it\Phi}(\lambda)$ satisfies the following formula:
 \begin{equation}\label{eqn2.01}
 {\it\Phi}(\lambda)=\mathtt{i}a\lambda+\frac{1}{2}\sigma^2\lambda^2+\int_{\mathbb{R}}(1-e^{\mathtt{i}\lambda u}+\mathtt{i}\lambda u)\nu(du),
 \end{equation}
 where $a\in\mathbb{R}$, $\sigma\geq 0$ and $\nu(du)$ is called \textit{L\'evy measure} with $$\int_{\mathbb{R}}(u\wedge u^2)\nu(du)<\infty.$$
 This assumption induces that the first moment of $\xi$ exists, i.e. $a=- \mathbf{E}[\xi_1]$.
 For any probability measure $\mu$ on $\mathbb{R}$, we denote by $\mathbf{P}_\mu$ and $\mathbf{E}_\mu$ the law and expectation of the L\'evy process $\xi$ started from $\mu$, respectively. When $\mu=\delta_x$ is a Dirac measure at point $x$, we write $\mathbf{P}_x$ for $ \mathbf{P}_{\delta_x}$ and $\mathbf{E}_x$ for $ \mathbf{E}_{\delta_x}$.
 For simplicity, we also write $\mathbf{P}$ for $ \mathbf{P}_0$ and $\mathbf{E}$ for $ \mathbf{E}_0$.

 We write $S=\{S_t:t\geq 0\}$ and $I=\{ I_t:t\geq 0 \}$ for the  \textit{running supremum} and \textit{infimum processes},
 \begin{equation*}
 S_t:= \sup\{0\vee \xi_s: 0\leq s\leq t\}
 \quad\mbox{and}\quad
  I_t:= \inf\{0\wedge \xi_s: 0\leq s\leq t\}.
 \end{equation*}
 Let $S-\xi := \{S_t-\xi_t: t\ge 0\}$ and $\xi-I := \{\xi_t-I_t: t\ge 0\}$ be the \textit{reflected processes} at the supremum and infimum respectively, which are Markov processes with Feller transition semigroups; see Proposition~1 in \cite{Bertoin1996}[p.156].
 For any $t>0$, the last passage times by $\xi$ at its supremum and infimum before $t$ are defined by
 \begin{equation*}
 G_t:=\sup\{s\leq t: \xi_s=S_t \mbox{ or }\xi_{s-} =S_t\}
  \quad\mbox{and}\quad
 g_t:=\sup\{s\leq t: \xi_s=I_t\mbox{ or }\xi_{s-}=I_t\}.
 \end{equation*}
 For any $x\leq 0$, we also define the first entrance times of $\xi$ in $(-\infty, x)$ and $(-x,\infty)$ by
 \begin{equation*}
 \tau_{x}:=\inf\{t> 0:\xi_t< x\}
 \quad\mbox{and}\quad  \tau_{-x}^+:=\inf\{t> 0:\xi_t>-x\}.
 \end{equation*}

 Let $L = \{L_t: t\ge 0\}$ be the local time of $S-\xi$ at zero in the sense of \cite{Bertoin1996}[p.109].
 Its \textit{inverse local time process} $L^{-1} = \{L^{-1}_t: t\ge 0\}$ is defined by
 \begin{equation*}
 L^{-1}_t :=\Big\{
 \begin{array}{ll}\inf\{s>0: L_s>t\}, & t< L_\infty;\cr
 \infty, & \mbox{otherwise}.
 \end{array}
 \end{equation*}
 The \textit{ladder height process} $H=\{H_t:t\ge 0\}$ of $\xi$ is defined by
 \begin{equation*}
 H_t :=\Big\{\begin{array}{ll}\xi_{L^{-1}_t}, & t< L_\infty;\cr
 \infty, & \mbox{otherwise}.
 \end{array}
 \end{equation*}
 By Lemma~2 in  \cite{Bertoin1996}[p.157], the two-dimensional process $(L^{-1},H)$ is a L\'{e}vy process (possibly killed at an exponential rate) and is well known as the \textit{ladder process} of $\xi$. It is usually characterized as follows:  for any $\lambda,u\ge 0$,
 \begin{equation}\label{eqn1.02}
 \mathbf{E}\big[\exp\{-\lambda L^{-1}_1 - u H_1\}\big]= \exp\{- \kappa(\lambda,u)\},
 \end{equation}
 where the \textit{bivariate exponent} $\kappa(\lambda,u)$ is given by
 \begin{equation}\label{eqn1.03}
 \kappa(\lambda,u)
  = k\exp\Big\{\int_0^\infty \frac{dt}{t}\int_{[0,\infty)}(e^{-t}-e^{-\lambda t-u x})\mathbf{P}\{\xi_t\in dx\}\Big\};
 \end{equation}
 see Corollary~10 in  \cite{Bertoin1996}[p.165-166].
 Here the constant $k>0$ is determined by the normalization of the local time.
 The \textit{renewal function} $V$ associated to the ladder height process $H$ is defined by
 \begin{equation}\label{eq2.6}
 V(x)
 = \int_0^\infty \mathbf{P}\{H_t\le x\}dt
 =\mathbf{E}\Big[\int_{[0,\infty)} 1_{\{S_t\le x\}}dL_t\Big], \qquad x\ge 0;
 \end{equation}
 see, \cite{Bertoin1996}[p.171].
 For the reflected process $\xi-I$, we can define the local time at $0$, the inverse local time process, the ladder height process and the renewal function in the same way as for $S-\xi$. They are denoted by $\hat L$, $\hat L^{-1}$, $\hat{H}$ and $\hat{V}$ respectively. Denote by $\hat \kappa(\lambda,u)$ the characteristic exponent of the ladder process $(\hat L^{-1},\hat H)$ with constant $\hat k>0$.  In this work, we choose some suitable normalization of the local times such that $k=\hat k=1$.

 From (\ref{eqn2.01}), the L\'evy process $\xi$ starting from $\xi_0$ admits the following L\'evy-It\^o's decomposition:
 \begin{equation}\label{eqn2.05}
 \xi_t=\xi_0 - at + \sigma B_t  +\int_0^t \int_{\mathbb{R}} u \tilde{N}(ds,du),
 \end{equation}
 where $B$ is a standard Brownian motion, $N(ds,du)$ is a Poisson random measure with intensity $ds\nu(du)$ and $\tilde{N}(ds,du):=N(ds,du)-ds\nu(du)$.
 For any $x>0$, we define $\xi^x =\{ \xi^x_t:t\geq 0  \}$ by removing from $\xi$ all jumps larger than $x$, i.e.
 \begin{equation}\label{eqn2.06}
 \xi^x_t:= \xi_0-\Big(a+\int_x^\infty u\nu(du)\Big)t +\sigma B_t+\int_0^t\int_{-\infty}^x u\tilde{N}(ds,du) ,
 \end{equation}
 which again is a L\'evy process with  characteristic exponent ${\it\Phi}^x(\lambda):=-\log\mathbf{E}[\exp\{\mathtt{i}\lambda\xi^x_1\}]$. Moreover, ${\it\Phi}^x(\lambda)$ satisfies  the following formula:
 \begin{equation}\label{eqn2.07}
 {\it\Phi}^x(\lambda)=\mathtt{i}\Big(a+\int_x^\infty u\nu(du)\Big)\lambda+\frac{1}{2}\sigma^2\lambda^2+\int_{-\infty}^x(1-e^{\mathtt{i}\lambda u}+\mathtt{i}\lambda u)\nu(du).
 \end{equation}


 \section{Asymptotic results for L\'evy processes with negative drift}
 \setcounter{equation}{0}

 In this section, under the following condition we provide several conditional limit theorems for L\'evy process $\xi$ with negative drift ($a>0$), which is equivalent that $\xi$ drifts to $-\infty$.
 \begin{condition}\label{Con2.1}
  There exit a constant $\alpha>1$ and a slowly varying function $\ell(x)$ at $\infty$ such that as $x\to \infty$,
  \begin{equation}\label{eqn3.01}
  \mathbf{P}\{\xi_1> x\} \sim x^{-\alpha} \cdot \ell(x).
  \end{equation}
 \end{condition}
 Actually, from (\ref{eqn2.05})-(\ref{eqn2.06}), we see that the right-tail of $\xi_1^1$ is light and hence $\xi_1\geq x$ for large $x$ if and only if there are large jumps happen in $[0,1]$, i.e.,
 \begin{equation*}
 \mathbf{P}\{\xi_1> x\} =\mathbf{P}\Big\{\xi^1_1+ \int_0^1\int_{1}^\infty u N(ds,du) > x\Big\} \sim \mathbf{P}\Big\{ \int_0^1\int_{1}^\infty u N(ds,du) > x\Big\}.
 \end{equation*}
 From Proposition 4.1 in \cite{Fay2006}, we have
 $\mathbf{P}\{\xi_1> x\} \sim 1-\exp\{-\nu(x,\infty)\}$ and hence Condition~\ref{Con2.1}  holds if and only if as $x\to\infty$,
 \begin{equation}\label{eqn3.02}
 \bar{\nu}(x):= \nu(x,\infty) \sim x^{-\alpha} \cdot \ell(x).
 \end{equation}

 \subsection{Limit theorems conditioned on $\xi_t\geq 0$}

 In this section, we study the number and the location of large jumps conditioned on $\xi_t>0$ for large $t$.
 For any $x,t>0$, let $\mathcal{N}_t^x$ be the number of jumps of $\xi$ larger than $x$ up to time $t$, i.e.,
 \begin{equation}\label{eqn3.03}
 \mathcal{N}_t^x :=\#\{s\leq t:\Delta\xi_s:=\xi_s-\xi_{s-}> x\} = \int_0^t \int_x^\infty N(ds,du).
 \end{equation}
 Let $\mathcal{J}^x :=\inf\{ s\geq 0:  \Delta\xi_s > x\}$ denote the arrival time of the first jump larger than $x$. Then for any $k\in\mathbb{N}$ and $t>0$,
 \begin{equation}\label{eqn3.04}
 \mathbf{P}\{\mathcal{N}_t^x>0\}=1-\exp\{- \bar\nu(x)\cdot t\}
 \quad \mbox{and}\quad
  \mathbf{P}\{\mathcal{N}_t^x=k\}= \frac{|\bar\nu(x)\cdot t|^k}{k !} \cdot \exp\{- \bar\nu(x)\cdot t\}.
 \end{equation}
 From these and Condition~\ref{Con2.1}, for any $b>0$ we have
 $\bar\nu(bt)\cdot t\to 0$ as $t\to\infty$ and hence
 \begin{equation}\label{eqn3.05}
 \mathbf{P}\{\mathcal{N}_t^{bt}>0\}\sim \mathbf{P}\{\mathcal{N}_t^{bt}=1\}\sim  \bar\nu(bt)\cdot t.
 \end{equation}
 Theorem~4 in \cite{DoneyJones2012} provided the following  asymptotic equivalences: as $t\to\infty$,
 \begin{equation}\label{eqn3.06}
 \mathbf{P}\{\xi_t> 0\}\sim t\cdot \mathbf{P}\{\xi_1> at\}\sim  \bar\nu(at)\cdot t.
 \end{equation}
 The following proposition proves that $\xi_t> 0$ for  large $t$ if and only if  there is one jump larger than $at$  in the time interval $[0,t]$.

 \begin{proposition}\label{Thm102}
 The two events $\mathcal{N}_t^{at}=1$ and $\xi_t> 0$ are asymptotically equivalent as $t\to\infty$, i.e.,
 \begin{equation}\label{eqn3.07}
 \mathbf{P}\big\{\mathcal{N}_t^{at}=1\,\big|\,\xi_t> 0\big\}\sim \mathbf{P}\big\{\xi_t> 0\,\big|\,\mathcal{N}_t^{at}=1\big\}\to 1.
 \end{equation}
 \end{proposition}
 \proof From Condition~\ref{Con2.1}, for any $\epsilon>0$, there exists $t_\epsilon>0$ such that $\int_{at_\epsilon}^\infty u\nu(du)\leq \epsilon$. From (\ref{eqn2.05}) and (\ref{eqn2.06}), for any $b\in(a+2\epsilon,a+3\epsilon)$ and $t>t_\epsilon$,
 conditioned on the event $\mathcal{N}_t^{bt}=1$ we have
 $ \xi_t\geq\xi^{bt_\epsilon}_t +\Delta \xi_{\mathcal{J}^{bt}} \geq \xi^{bt_\epsilon}_t +bt.$
 From (\ref{eqn2.06}), (\ref{eqn3.03}) and the spatial orthogonality of $N(ds,du)$, for any $t> t_\epsilon$ we have $\xi^{b t_\epsilon}_t $ is independent of $\mathcal{N}_t^{bt}$  and
 \begin{equation}\label{eqn3.07.01}
  \mathbf{P}\big\{\xi_t> 0, \mathcal{N}_t^{at}=1\big\}\geq \mathbf{P}\big\{\xi^{b t_\epsilon}_t+bt> 0, \mathcal{N}_t^{bt}=1\big\}
 = \mathbf{P}\big\{\xi^{b t_\epsilon}_t+bt> 0\big\}\cdot\mathbf{P}\big\{\mathcal{N}_t^{bt}=1\big\}.
 \end{equation}
 Notice that $\mathbf{E}\big[\xi^{b t_\epsilon}_1+b\big] \geq \epsilon$ and hence  $\mathbf{P}\big\{\xi^{b t_\epsilon}_t+bt> 0\big\}\to 1$ as $t\to\infty$.
 From (\ref{eqn3.05}),
 \begin{eqnarray*}
  \liminf_{t\to\infty}    \mathbf{P}\big\{\xi_t> 0\,\big|\,  \mathcal{N}_t^{at}=1\big\}
  \ar\geq\ar\lim_{b\to a+} \liminf_{t\to\infty} \frac{ \mathbf{P}\big\{\xi_t> 0, \mathcal{N}_t^{at}=1\big\}}{\mathbf{P}\big\{  \mathcal{N}_t^{bt}=1\big\}}\frac{\mathbf{P}\big\{  \mathcal{N}_t^{bt}=1\big\}}{ \mathbf{P}\big\{  \mathcal{N}_t^{at}=1\big\}}\geq \lim_{b\to a+}\Big| \frac{a}{b}\Big|^{\alpha}=1,
 \end{eqnarray*}
which induces that $ \mathbf{P}\big\{\xi_t> 0| \mathcal{N}_t^{at}=1\big\}\to 1$ as $t\to\infty$.
 On the other hand, from (\ref{eqn3.05}), (\ref{eqn3.06})  and (\ref{eqn3.07.01}) we also have
  \begin{eqnarray*}
 	\liminf_{t\to\infty}    \mathbf{P}\big\{ \mathcal{N}_t^{at}=1\,\big|\, \xi_t> 0 \big\}
 	\ar\geq\ar\lim_{b\to a+} \liminf_{t\to\infty} \frac{ \mathbf{P}\big\{\xi_t> 0, \mathcal{N}_t^{at}=1\big\}}{\mathbf{P}\big\{  \mathcal{N}_t^{bt}=1\big\}}\frac{\mathbf{P}\big\{  \mathcal{N}_t^{bt}=1\big\}}{ \mathbf{P}\big\{  \xi_t> 0 \big\}}\geq \lim_{b\to a+}\Big| \frac{a}{b}\Big|^{\alpha}=1,
 \end{eqnarray*}
which immediately induces that $ \mathbf{P}\big\{ \mathcal{N}_t^{at}=1\,\big|\, \xi_t> 0 \big\}\to1$ as $t\to\infty$.
 Here we have finished the proof.
 \qed

 The following theorem shows that conditioned on $\xi_t> 0$ for large $t$,  the only jump larger than $at$ occurs at a random time that is uniformly distributed in $[0,t]$. This also can be gotten from Theorem~6 in \cite{DoneyJones2012}.
 \begin{theorem}\label{Thm103}
 There exists a uniformly distributed random variable $\mathbf{U}$ on $[0,1]$ such that $\big\{\mathcal{J}^{at}/t\,\big|\,\xi_t> 0\big\} \to\mathbf{U}$ in distribution as $t\to\infty$.
 \end{theorem}
 \proof For any $U>0$, it suffices to prove that $\mathbf{P}\big\{\mathcal{J}^{at}/t> U\,\big|\,\xi_t> 0\big\}\to (1-U)^+$. From Proposition~\ref{Thm102}, we have  for $t>0$ large enough,
 \begin{equation*}
 \mathbf{P}\big\{\mathcal{J}^{at}/t> U\,\big|\,\xi_t> 0\big\}\sim  \mathbf{P}\big\{\mathcal{J}^{at}/t> U\,\big|\,\mathcal{N}_t^{at}=1\big\} \sim \mathbf{P}\big\{\mathcal{N}^{at}_{Ut}=0\,\big|\,\mathcal{N}_t^{at}=1\big\},
 \end{equation*}
 which equals $0$ if $U>1$. From the independent increments of $\xi$ and (\ref{eqn3.04}), we have for $U\in [0,1]$,
 \begin{equation*}
 \mathbf{P}\big\{\mathcal{N}^{at}_{Ut}=0\,\big|\,\mathcal{N}_t^{at}=1\big\}
 = \frac{\mathbf{P}\big\{\mathcal{N}_{Ut}^{at}=0\big\}\cdot\mathbf{P}\big\{\mathcal{N}_{(1-U)t}^{at}=1 \big\}}{\mathbf{P}\big\{\mathcal{N}_t^{at}=1\big\}}
 =1-U.
 \end{equation*}
 \qed

 \subsection{Limit theorems conditioned to stay positive}

 In this section, we provide several conditional limit theorems for the L\'evy process $\xi$ conditioned to stay positive.
 For any $x>0$, from Proposition~17 in  \cite{Bertoin1996}[p.172]  we have $\mathbf{E}_x[\tau_{0}]= C\cdot \hat{V}(x)<\infty$ for some constant $C>0$.
 \cite{DenisovShneer2013} proved in Theorem~2.2 that as  $t\rightarrow \infty$,
 \begin{equation}\label{eqn3.09}
 \mathbf{P}_x\{\tau_0>t\}\sim \mathbf{E}_x[\tau_0]\cdot\frac{\mathbf{P}\{\xi_t> 0\}}{t} \sim \mathbf{E}_x[\tau_0] \cdot \mathbf{P}\{\xi_1> at\}\sim \mathbf{E}_x[\tau_0] \cdot \bar{\nu}(at).
 \end{equation}
 From Proposition~\ref{Thm102} and Theorem~\ref{Thm103}, we see that $\mathbf{P}_x\{ \xi_t> 0  \}\sim  \sum_{k=0}^{[t]} \mathbf{P}\{ \mathcal{J}^{at}\in(k,k+1] ,\xi_t> 0 \}  $ for large $t$ and $\mathbf{P}\{ \mathcal{J}^{at}\in(i,i+1],\xi_t> 0   \}\sim\mathbf{P}\{ \mathcal{J}^{at}\in(j,j+1] ,\xi_t> 0  \}$ for any $i,j\geq 0$.
 Comparing this to (\ref{eqn3.09}), we may conjecture that conditioned to stay positive the probability of $ \mathcal{J}^{at}\in(k,k+1]$ should decrease as $k\to\infty$.
 The following theorem shows that conditionally on staying positive the arrival time of the first large jump is distributed as a size-biased distribution, which induces that  the early arrival of the first large jump is necessary.

 \begin{theorem}\label{Thm106}
 For any $x>0$,  $T\geq 0$  and $b\geq a$, we have as $t\rightarrow\infty$,
 \begin{equation}\label{eqn3.11}
 \mathbf{P}_x\big\{\Delta\xi_{\mathcal{J}^{at}}>bt, \mathcal{J}^{at} \leq T\,|\,\tau_0>t\big\}\to (b/a)^{-\alpha}\cdot \frac{\mathbf{E}_x[\tau_0\wedge T]}{\mathbf{E}_x[\tau_0]}.
 \end{equation}
 \end{theorem}
 \proof We first prove this result with $b>a$.  Notice that $\mathcal{J}^{at}$ is a stopping time for any $t>0$. Let $\xi'$ be an independent copy of $\xi$.
 From the strong Markov property of $\xi$,
 \begin{eqnarray*}
 \lefteqn{ \mathbf{P}_x\big\{\Delta\xi_{\mathcal{J}^{at}}>bt, \mathcal{J}^{at}\leq T,\tau_0>t\big\}}\ar\ar\cr
 \ar\ar\cr
 \ar\ar=  \mathbf{P}_x\Big\{\Delta\xi_{\mathcal{J}^{at}}>bt, \mathcal{J}^{at}\leq T,\tau_0>\mathcal{J}^{at},\mathbf{P}_{\xi_{\mathcal{J}^{at}}}\Big\{\inf_{r\in[0,t-\mathcal{J}^{at}]}\xi'_r> 0\Big\} \Big\}.
 \end{eqnarray*}
 Conditioned on $\mathcal{J}^{at}\leq T$, $\tau_0>\mathcal{J}^{at}$ and $\Delta\xi_{\mathcal{J}^{at}}>bt$, for large $t$ we have $\xi_{\mathcal{J}^{at}}> bt$ and
 \begin{equation*}
 \mathbf{P}_{\xi_{\mathcal{J}^{at}}}\Big\{\inf_{r\in[0,t-\mathcal{J}^{at}]}\xi'_r\geq 0\Big\}
 \geq  \mathbf{P} \Big\{\inf_{r\in[0,t]}\xi'_r\geq -\xi_{\mathcal{J}^{at}}\Big\}
 \geq \mathbf{P}\{I_t> -bt \},
 \end{equation*}
 which goes to $1$ as $t\to\infty$ because of $\mathbf{E}[\xi_1]=a<b$ and hence
 \begin{eqnarray*}
 \mathbf{P}_x\big\{\Delta\xi_{\mathcal{J}^{at}}>bt, \mathcal{J}^{at}\leq T,\tau_0>t\big\}
 \ar\geq\ar \mathbf{P}\big\{I_t> -bt \big\} \cdot  \mathbf{P}_x\big\{\Delta\xi_{\mathcal{J}^{at}}>bt, \mathcal{J}^{at}\leq T,\tau_0>\mathcal{J}^{at}\big\} \cr
 \ar\ar\cr
 \ar\sim \ar \mathbf{P}_x\big\{\Delta\xi_{\mathcal{J}^{at}}>bt, \mathcal{J}^{at}\leq T,\tau_0>\mathcal{J}^{at}\big\}.
 \end{eqnarray*}
 Together with the fact that
 $	\mathbf{P}_x\big\{\Delta\xi_{\mathcal{J}^{at}}>bt, \mathcal{J}^{at}\leq T,\tau_0>t\big\}\leq  \mathbf{P}_x\big\{\Delta\xi_{\mathcal{J}^{at}}>bt, \mathcal{J}^{at}\leq T,\tau_0>\mathcal{J}^{at}\big\}$
 for any $t> T$, we have as $t\to\infty$,
 \begin{equation*}
 \mathbf{P}_x\big\{\Delta\xi_{\mathcal{J}^{at}}>bt, \mathcal{J}^{at}\leq T,\tau_0>t\big\}\sim \mathbf{P}_x\big\{\Delta\xi_{\mathcal{J}^{at}}>bt, \mathcal{J}^{at}\leq T,\tau_0>\mathcal{J}^{at}\big\}.
 \end{equation*}
 To get the desired result, it suffices to prove that
 \begin{equation*}
 \lim_{t\rightarrow\infty}\frac{ \mathbf{P}_x\big\{\Delta\xi_{\mathcal{J}^{at}}>bt, \mathcal{J}^{at}\leq T,\tau_0>\mathcal{J}^{at}\big\}}{\mathbf{P}_x\{\tau_0>t\}}
 =(b/a)^{-\alpha} \cdot \frac{\mathbf{E}_x[\tau_0 \wedge T]}{\mathbf{E}_x[\tau_0]}.
 \end{equation*}
 From (\ref{eqn3.04}) and (\ref{eqn3.09}), we first have as  $t\to\infty$,
 \begin{equation*}
 \mathbf{P}_x\{\mathcal{N}_T^{bt}>1\}
 =1- e^{-\bar{\nu}(bt)\cdot T}- \bar\nu(bt)T e^{- \bar\nu(bt)\cdot T} \sim \frac{|\bar\nu(bt)\cdot T|^2}{2} =o\Big(\mathbf{P}_x\{\tau_0>t\}\Big).
 \end{equation*}
 Moreover, we also have
\begin{equation*}
 \big\{\Delta\xi_{\mathcal{J}^{at}}>bt, \mathcal{J}^{at}\leq T,\tau_0>\mathcal{J}^{at}\big\}
 = \big\{ \mathcal{J}^{at}= \mathcal{J}^{bt}, \mathcal{J}^{at}\leq T,\tau_0>\mathcal{J}^{at}\big\}
 \subset
 \big\{ \mathcal{N}_T^{bt}\geq 1,  \tau_0>\mathcal{J}^{at} \big\}
 \end{equation*}
 and
 \begin{equation*}
 \big\{\Delta\xi_{\mathcal{J}^{at}}>bt, \mathcal{J}^{at}\leq T,\tau_0>\mathcal{J}^{at}\}\cup\{ \mathcal{N}_T^{at}>1 \big\}
 \supset
 \big\{ \mathcal{N}_T^{bt}= 1,  \tau_0>\mathcal{J}^{at} \big\}.
 \end{equation*}
 These immediately induce that
 \begin{equation}\label{eqn3.12}
 \frac{ \mathbf{P}_x\{\Delta\xi_{\mathcal{J}^{at}}>bt, \mathcal{J}^{at}\leq T,\tau_0>\mathcal{J}^{at}\}}{\mathbf{P}_x\{\tau_0>t\}}
 \sim
 \frac{ \mathbf{P}_x\{\mathcal{N}_T^{bt}\geq 1,\tau_0>\mathcal{J}^{bt}\}}{\mathbf{P}_x\{\tau_0>t\}}.
 \end{equation}
 Let $I^{bt}=\{ I^{bt}_s:s\geq 0 \}$ be the running infimum process of $\xi^{bt}$.
 For any fixed $t$, we see that $\xi^{bt}$ is independent of $\mathcal{J}^{bt}$ and
 $\xi_s^{bt}=\xi_s$ for any $s<\mathcal{J}^{bt}$. These induce that $\tau_0>\mathcal{J}^{bt} $ if and only if $I_s=I_s^{bt}>0$ for any $s\leq \mathcal{J}^{bt}$.
 Thus
 \begin{eqnarray*}
 	\big\{\mathcal{N}_T^{bt}\geq 1,\tau_0>\mathcal{J}^{bt}\big\}
 	=\Big\{  \int_0^T \int_{bt} \mathbf{1}_{\{ I_r^{bt}\geq 0\} } N(dr,dx)\geq 1 \Big\}.
 	\end{eqnarray*}
 From the independence between $I_\cdot^{bt}$ and $\int_0^\cdot \int_{bt}^\infty N(dr,dx) $ for any $t>0$, we have
 \begin{eqnarray*}
  \int_0^T \int_{bt} \mathbf{1}_{\{ I_r^{bt}\geq 0\} } N(dr,dx)= \mathbf{N}_T
 \end{eqnarray*}
in distribution, where $ \{\mathbf{N}_r:r\geq 0\}$ is a general Poisson point process with conditional intensity $\bar{\nu}(bt) \mathbf{1}_{\{ I_r^{bt}\geq 0\} }$.
 From this and the properties of general Poisson point process, we have as $t\to\infty$,
 \begin{eqnarray}\label{eqn3.13}
 \mathbf{P}_x\big\{\mathcal{N}_T^{bt}\geq 1,\tau_0>\mathcal{J}^{bt}\big\} = \mathbf{P}_x\{ \mathbf{N}_T\geq 1 \}
 \ar=\ar 1-\mathbf{E}_x\Big[ \exp\Big\{- \bar{\nu}(bt) \int_0^T\mathbf{1}_{\{ I_r^{bt}\geq 0\} }  dr   \Big\} \Big] \cr
 \ar\sim\ar  \bar{\nu}(bt) \int_0^T \mathbf{P}_x\big\{ I_r^{bt}\geq 0\big\} dr .
 \end{eqnarray}
 Since $\sup_{s\in[0,T]}|I_s^{bt}-I_s|\to 0$ a.s. as $t\to\infty$, we have
 \begin{equation}\label{eqn3.14}
 \mathbf{P}_x\big\{\mathcal{N}_T^{bt}\geq 1,\tau_0>\mathcal{J}^{bt}\big\}
 \sim \bar{\nu}(bt) \int_0^T \mathbf{P}_x\{ I_r\geq 0\} dr
 = \bar{\nu}(bt)\cdot \int_0^T \mathbf{P}_x\{\tau_0> r\}dr.
 \end{equation}
 Taking this back into (\ref{eqn3.12}), from (\ref{eqn3.02}) and (\ref{eqn3.09}) we have as $t\to\infty$,
 \begin{eqnarray}\label{eqn3.15}
  \frac{ \mathbf{P}_x\{\Delta\xi_{\mathcal{J}^{at}}>bt, \mathcal{J}^{at}\leq T,\tau_0>\mathcal{J}^{at}\}}{\mathbf{P}_x\{\tau_0>t\}}
 \ar\sim\ar
 (b/a)^{-\alpha} \cdot \frac{\mathbf{E}_x[\tau_0 \wedge T]}{\mathbf{E}_x[\tau_0]}.
 \end{eqnarray}
 Now we consider the case with $b=a$. From the previous result, we have
 \begin{eqnarray*}
\liminf_{t\to\infty} \mathbf{P}_x\{\mathcal{J}^{at} \leq T\,\big|\,\tau_0>t\}\ar\geq\ar  \lim_{c\to a+}\lim_{t\to\infty}\mathbf{P}_x\{\Delta\xi_{\mathcal{J}^{at}}>ct, \mathcal{J}^{at} \leq T \,\big|\, \tau_0>t\}=  \frac{\mathbf{E}_x[\tau_0 \wedge  T]}{\mathbf{E}_x[\tau_0]}.
 \end{eqnarray*}
 Moreover,  like the deduction in (\ref{eqn3.13})-(\ref{eqn3.15}) we also have as $t\to\infty$
 \begin{eqnarray*}
 \limsup_{t\to\infty} \mathbf{P}_x\big\{\mathcal{J}^{at} \leq T \,\big|\, \tau_0>t\big\}
 \ar\leq\ar \lim_{c\to a-}\limsup_{t\to\infty} \mathbf{P}_x\big\{\mathcal{J}^{ct} \leq T \,\big|\, \tau_0>t\big\} \cr
 \ar\leq\ar  \lim_{c\to a-}\limsup_{t\to\infty}\frac{\mathbf{P}_x\{\mathcal{N}_T^{ct}\geq 1,\tau_0>\mathcal{J}^{ct}\}}{\mathbf{P}_x\{ \tau_0>t\}}\cr
 \ar=\ar \lim_{c\to a-}  (c/a)^{-\alpha} \cdot \frac{\mathbf{E}_x[\tau_0\wedge  T]}{\mathbf{E}_x[\tau_0]}=\frac{\mathbf{E}_x[\tau_0\wedge T]}{\mathbf{E}_x[\tau_0]}.
 \end{eqnarray*}
 The desired result for $b=a$ follows directly from these two results.
 \qed

 For any large $t$, from (\ref{eqn3.09}) we see that $\mathbf{P}_x\{\tau_0>t\}\to1$ as $x\to \infty$, which means that $\xi$ will stay positive for a long time after jumping into a set far away from the origin.
 This recommends us that  $\xi$ would stay positive for a long time if and only if it can stay  positive until the arrival of the first large jump; see the following lemma.
 \begin{lemma}\label{Thm107}
 For any $x> 0$ and large $t$, we have
  as  $t\to\infty$,
 \begin{equation}\label{eqn3.16}
 \mathbf{P}_x\big\{\tau_0> \mathcal{J}^{at}\,\big|\,\tau_0> t\big\}
 \sim \mathbf{P}_x\big\{\tau_0> t\,\big|\,\tau_0> \mathcal{J}^{at}\big\} \to 1.
 \end{equation}
 \end{lemma}
 \proof Obviously, it suffices to prove that
 $ \mathbf{P}_x\big\{\{\tau_0> \mathcal{J}^{at}\}\Delta \{\tau_0> t\}\big\}=o\big(\mathbf{P}_x\{\tau_0> t\}\big).$
 From the properties of symmetric difference, we have
 \begin{eqnarray*}
 \mathbf{P}_x\big\{\{\tau_0> \mathcal{J}^{at}\}\Delta \{\tau_0> t\}\big\}
 \ar=\ar \mathbf{P}_x\big\{ \tau_0> t\big\}-\mathbf{P}_x\big\{\tau_0> \mathcal{J}^{at}\big\}
+2 \mathbf{P}_x\big\{ \mathcal{J}^{at}<\tau_0\leq t\big\}.
 \end{eqnarray*}
  Let $I^{at}=\{ I^{at}_s:s\geq 0 \}$ be the running infimum process of $\xi^{at}$.
  Like the computations in (\ref{eqn3.13})-(\ref{eqn3.14}), from the fact that $\mathbf{E}_x[\int_0^\infty \mathbf{1}_{\{ I_s\geq 0 \}}ds]=\mathbf{E}_x[\tau_0]<\infty$
  we have as $t\to\infty$,
   \begin{eqnarray*}
  	  \mathbf{P}_x\big\{\tau_0>\mathcal{J}^{at}\big\}
  	  \ar=\ar \mathbf{P}  \Big\{ \int_0^\infty \int_{at}^\infty \mathbf{1}_{\{I^{at}_r>0\}}N(dr,dx)\geq 0  \Big\}\cr
  	  \ar=\ar 1-\mathbf{E}_x\Big[ \exp\Big\{- \bar{\nu}(bt) \int_0^\infty \mathbf{1}_{\{ I_r^{bt}>0\} }  dr   \Big\} \Big] \cr
  	  \ar\leq\ar  1-\mathbf{E}_x\Big[ \exp\Big\{- \bar{\nu}(bt) \int_0^\infty \mathbf{1}_{\{ I_r>0\} }  dr   \Big\} \Big]
  	  \sim \mathbf{E}_x[\tau_0]\cdot \bar{\nu}(bt).
  	 \end{eqnarray*}
   On the other side, we also have
  \begin{eqnarray*}
  \mathbf{P}_x\big\{\tau_0>\mathcal{J}^{at}\big\} \geq 	\mathbf{P}_x\{\tau_0>\mathcal{J}^{at}, \mathcal{J}^{at}\leq t\}
 	\ar=\ar\mathbf{P}_x\{\tau_0>\mathcal{J}^{at},\, \mathcal{N}^{at}_t\geq 1\}\cr
 	\ar\sim \ar \bar\nu(at)\int_0^\infty\mathbf{P}_x\{\tau_0> s\}ds
 	\sim \mathbf{E}_x[\tau_0] \cdot \bar\nu(at).
 \end{eqnarray*}
Putting these two estimates together with (\ref{eqn3.09}), we have
 \begin{equation}\label{eqn3.17}
 \lim_{t\to \infty}
\frac{\mathbf{P}_x\big\{\tau_0>\mathcal{J}^{at}\big\}}{\mathbf{P}_x\{\tau_0>t\}}
= \lim_{t\to\infty}\frac{\mathbf{P}_x\{\tau_0>\mathcal{J}^{at}, \,\mathcal{J}^{at}\leq t\}}{\mathbf{P}_x\{\tau_0>t\}}  =1.
\end{equation}
Moreover,  applying Theorem~\ref{Thm106} with $b=a$, for any $T>0$ we have
\begin{eqnarray*}
	\liminf_{t\to\infty}\frac{\mathbf{P}_x\{\mathcal{J}^{at}\leq t< \tau_0\}}{\mathbf{P}_x\{\tau_0>t\}}
	\ar\geq\ar  \lim_{t\rightarrow\infty}\mathbf{P}_x\{ \mathcal{J}^{at} \leq T \,\big|\, \tau_0>t\}
	=  \frac{\mathbf{E}_x[\tau_0 \wedge T]}{\mathbf{E}_x[\tau_0]},	
\end{eqnarray*}
which goes to $1$ as $T\to\infty$.
This result together with (\ref{eqn3.17}) immediately induces that
\begin{equation*}
\lim_{t\to \infty}\frac{\mathbf{P}_x\{ \mathcal{J}^{at}<\tau_0\leq t\}}{\mathbf{P}_x\{\tau_0>t\}}
= \lim_{t\to \infty}\frac{\mathbf{P}_x\{\tau_0> \mathcal{J}^{at},\,\mathcal{J}^{at}\leq t\}}{\mathbf{P}_x\{\tau_0>t\}}
-\lim_{t\to\infty}\frac{\mathbf{P}_x\{\mathcal{J}^{at}\leq t< \tau_0\}}{\mathbf{P}_x\{\tau_0>t\}}
=0.
\end{equation*}
Combining this result with (\ref{eqn3.17}), we can immediately get the desired result, i.e.,
\begin{equation*}
\frac{\mathbf{P}_x\left\{\{\tau_0> \mathcal{J}^{at}\}\Delta \{\tau_0> t\}\right\}}{\mathbf{P}_x\{\tau_0> t\}}
= 1-\frac{\mathbf{P}_x\big\{\tau_0> \mathcal{J}^{at}\big\} }{\mathbf{P}_x\big\{ \tau_0> t\big\}}
+ \frac{2\mathbf{P}_x\big\{ \mathcal{J}^{at}<\tau_0\leq t\big\}}{\mathbf{P}_x\big\{ \tau_0> t\big\}},
\end{equation*}
which vanishes as $t\to\infty$.
 \qed

 From the previous results, we see that if $\xi$ could stay positive up to a large time $t$, there must have been a jump larger than $at$ which has occurred very early.
 Moreover, after the large jumps the sample paths will stay at the high position for a long time.
 To describe this phenomena,  the following theorem provides a limit theorem for the spatial-scaled process conditioned to stay positive.
 It shows that compared to the large jump, the effect of the downward drift on the process in the future can almost be ignored.
 The similar but different discrete version of this theorem for random walks with finite variance was established by Theorem~3.2 in \cite{Durrett1980}.
  \begin{theorem}\label{Thm108}
  	Fix $x>0$, let $\mathcal{T}_x$ and $\mathcal{P}_{a,\alpha}$ be two independent positive random variables with $\mathbf{P}\{\mathcal {T}_x\leq t\}=\mathbf{E}_x[\tau_0\wedge  t]/\mathbf{E}_x[\tau_0]$ for any $t\geq0$
  	and
  	$\mathbf{P}\{\mathcal{P}_{a,\alpha}\geq z\}=(z/a)^{-\alpha}$ for any $z\geq a$.
  	Then $\{t^{-1}\xi_s,s\geq 0|\tau_0>t, \xi_0=x\}$ converges to $\{\mathcal{P}_{a,\alpha}\cdot \mathbf{1}_{\{\mathcal {T}_x<s\}},s\geq 0\}$ weakly  in $\mathbf{D}([0,\infty),\mathbb{R})$ as $t\to\infty$.
  \end{theorem}
 \proof  It is easy to see that the desired result follows directly from the following two statements:
 \begin{enumerate}	
 \item[(i)] For any fixed $T>0$ and $\varepsilon>0$, we have as $t\to\infty$,
   \begin{equation*}
   \lim_{t\to\infty}\mathbf{P}_x\Big\{\sup_{s\in[0,T]}\big|t^{-1}\xi(s) - t^{-1}\Delta \xi_{\mathcal{J}^{at}}\cdot \mathbf{1}_{[\mathcal{J}^{at},\infty)}(s) \big|\geq \varepsilon \  \Big|\ \tau_0> t\Big\}\to 0;
   \end{equation*}

 \item[(ii)] As $t\rightarrow\infty$ we have $\big\{t^{-1}\Delta \xi_{\mathcal{J}^{at}}\cdot \mathbf{1}_{[\mathcal{J}^{at},\infty)}(s):s\geq 0 \,\big|\,\tau_0>t,\xi_0=x \big\}$ converges weakly to $\big\{\mathcal{P}_{a,\alpha}\cdot \mathbf{1}_{[\mathcal{T}_x,\infty)}(s):s\geq 0\big\}$  in $\mathbf{D}([0,\infty),\mathbb{R})$.

 \end{enumerate}
 For (i), we first decompose $\xi $ at the stopping time $\mathcal{J}^{at}$ as follows: for any $s\geq 0$,
 \begin{equation*}
 \xi_s=\xi_{s\wedge \mathcal{J}^{at}-}+\Delta\xi_{\mathcal{J}^{at}}\mathbf{1}_{[\mathcal{J}^{at},\infty)}(s)+ \xi_{s\vee \mathcal{J}^{at}}-\xi_{\mathcal{J}^{at}}.
 \end{equation*}
 From this and Lemma~\ref{Thm107},  we have for $t>0$ large enough,
 \begin{eqnarray}\label{eqn3.20}
  \lefteqn{\mathbf{P}_x\Big\{\sup_{s\in[0,T]}\big| \xi(s)- \Delta \xi_{\mathcal{J}^{at}}\mathbf{1}_{[\mathcal{J}^{at},\infty)}(s)\big| \geq \varepsilon t \,\Big|\, \tau_0> t\Big\}}\ar\ar\cr
  \ar\sim \ar \mathbf{P}_x\Big\{\sup_{s\in[0,T]} \big|\xi_{s\wedge \mathcal{J}^{at}-}+ \xi_{s\vee \mathcal{J}^{at}}-\xi_{\mathcal{J}^{at}}\big| \geq \varepsilon t \,\Big|\, \tau_0>\mathcal{J}^{at} \Big\} \cr
  \ar\leq\ar   \mathbf{P}_x\Big\{\sup_{s\in[0,T]} |  \xi_{s\vee \mathcal{J}^{at}}-\xi_{\mathcal{J}^{at}}| \geq \varepsilon t/2\mbox{\quad or \ }  \sup_{s\in[0,T]} |\xi_{s\wedge J^{at}-}| \geq \varepsilon t/2 \,\Big|\, \tau_0>\mathcal{J}^{at}\Big\}.
 \end{eqnarray}
 By the strong Markov property and the independent increments of  $\xi$,
 \begin{eqnarray*}
 \mathbf{P}_x\Big\{\sup_{s\in[0,T]} |  \xi_{s\vee \mathcal{J}^{at}}-\xi_{\mathcal{J}^{at}}| \geq \varepsilon t/2 \,\Big|\, \tau_0>\mathcal{J}^{at}\Big\}
 \ar=\ar \mathbf{P}\Big\{\sup_{s\in[0,T]} |  \xi_{s\vee \mathcal{J}^{at}}-\xi_{\mathcal{J}^{at}}| \geq \varepsilon t/2   \Big\}\cr
 \ar\leq\ar \mathbf{P}\Big\{\sup_{s\in[0,T]} |\xi_s|>\varepsilon t/2\Big\},
 \end{eqnarray*}
 which vanishes as $t\to\infty$.
 From the fact that $\xi_{s}=\xi^{at}_{s}$ for any $s<\mathcal{J}^{at}$ and the independence between $\xi^{at}$ and $\mathcal{J}^{at}$, we have for any $M>0$,
 \begin{eqnarray*}
  \lefteqn{\mathbf{P}_x\Big\{\sup_{s\in[0,T]} |\xi_{s\wedge J^{at}-}| \geq \varepsilon t/2,  \tau_0>\mathcal{J}^{at}\Big\}}\ar\ar\cr
  \ar\leq\ar   \mathbf{P}_x\Big\{\sup_{s\in[0,T]} |\xi_{s\wedge J^{at}-}| \geq \varepsilon t/2,  \mathcal{J}^{at}<M\Big\} +  \mathbf{P}_x\big\{  \mathcal{J}^{at}>M, \tau_0>\mathcal{J}^{at}\big\}\cr
  \ar\leq\ar \mathbf{P}_x\Big\{\sup_{s\in[0,T]} |\xi_{s}^{at}| \geq \varepsilon t/2\Big\}\cdot \mathbf{P}\{\mathcal{J}^{at}<M\} + \mathbf{P}_x\big\{  \mathcal{J}^{at}>M, \tau_0>\mathcal{J}^{at}\big\}.
 \end{eqnarray*}
 From (\ref{eqn3.04}), (\ref{eqn3.09}) and Lemma~\ref{Thm107},  we have $ \mathbf{P}\{\mathcal{J}^{at}<M\}=  \mathbf{P}\{\mathcal{N}_M^{at}\geq 1\}\sim M\bar\nu(at)\sim M \mathbf{P}_x\{\tau_0>\mathcal{J}^{at}\}/\mathbf{E}_x[\tau_0]$ as $t\to\infty$ and hence
 \begin{equation*}
 \mathbf{P}_x\Big\{\sup_{s\in[0,T]} |\xi_{s}^{at}| \geq \varepsilon t/2\Big\}\cdot \mathbf{P}\{\mathcal{J}^{at}<M\}  =o(\mathbf{P}\{\tau_0>\mathcal{J}^{at}\}).
 \end{equation*}
 From the property of the symmetric difference of two sets, we have
  \begin{eqnarray*}
\lefteqn{\mathbf{P}_x\big\{ \{ \mathcal{J}^{at}>M\} \cap\big(\{ \tau_0>\mathcal{J}^{at}\}\Delta \{\tau_0>t\}\big)\big\}}\ar\ar\cr
\ar\ar\cr
 \ar=\ar
 \mathbf{P}_x\big\{  \mathcal{J}^{at}>M, \tau_0>\mathcal{J}^{at}\big\}
 -\mathbf{P}_x\big\{  \mathcal{J}^{at}>M, \tau_0>t\big\}
+2\mathbf{P}_x\big\{  \mathcal{J}^{at}>M, t< \tau_0\leq \mathcal{J}^{at}\big\}.
 \end{eqnarray*}
 From Lemma~\ref{Thm107}, as $t\to\infty$ we have $\mathbf{P}_x\{ \tau_0>t \}\sim \mathbf{P}_x\big\{   \tau_0>\mathcal{J}^{at}\big\}$ and
  \begin{eqnarray*}
  \frac{\mathbf{P}_x\big\{ \{ \mathcal{J}^{at}>M\} \cap\big(\{ \tau_0>\mathcal{J}^{at}\}\Delta \{\tau_0>t\}\big)\big\}}{\mathbf{P}_x\{ \tau_0>t \}}
  	\ar\leq \ar   \frac{\mathbf{P}_x\big\{ \{   \tau_0>\mathcal{J}^{at}\}\Delta \{\tau_0>t\}\big\}}{\mathbf{P}_x\big\{   \tau_0>\mathcal{J}^{at}\big\}}
  	\to0.
  \end{eqnarray*}
 Moreover, from (\ref{eqn3.16}) we also have  as $t\to\infty$,
  \begin{eqnarray*}
 	\frac{\mathbf{P}_x\big\{  \mathcal{J}^{at}>M, t< \tau_0\leq \mathcal{J}^{at}\big\}}{\mathbf{P}_x\{ \tau_0>t \}}
 	\ar=\ar   \frac{\mathbf{P}_x\big\{   t< \tau_0\leq \mathcal{J}^{at}\big\}}{\mathbf{P}_x\{ \tau_0>t \}}
 	\to0.
 \end{eqnarray*}
 Putting these three results above together with Theorem~\ref{Thm106}, we have
  \begin{eqnarray*}
 \lim_{M\to\infty} \limsup_{t\to\infty}\frac{\mathbf{P}_x\big\{  \mathcal{J}^{at}>M, \tau_0>\mathcal{J}^{at}\big\}}{\mathbf{P}_x\big\{   \tau_0>\mathcal{J}^{at}\big\}}
 \ar=\ar    \lim_{M\to\infty}\limsup_{t\to\infty}\frac{\mathbf{P}_x\big\{  \mathcal{J}^{at}>M, \tau_0>t\big\}}{\mathbf{P}_x\big\{   \tau_0>t\big\}}\cr
 \ar\ar\cr
 \ar= \ar  \lim_{M\to\infty} \limsup_{t\to\infty} \mathbf{P}_x\big\{ \mathcal{J}^{at}>M \,\big|\,  \tau_0>t\big\}=0.
 \end{eqnarray*}
 Putting these estimates together, we have
 \begin{equation*}
 \lim_{t\to\infty} \mathbf{P}_x\Big\{\sup_{s\in[0,T]} |\xi_{s\wedge J^{at}-}| \geq \varepsilon t/2 \,\Big|\, \tau_0>\mathcal{J}^{at}\Big\}=0
 \end{equation*}
 and  (i) follows. Now we start to prove (ii).
 From Lemma~\ref{Thm107}, it suffices to prove that $ \big\{Y^t_s:=t^{-1}\Delta \xi_{\mathcal{J}^{at}}\cdot \mathbf{1}_{[\mathcal{J}^{at},\infty)}(s):s\in[0,1] \big|\tau_0>\mathcal{J}^{at} ,\xi_0=x\big\}$ converges weakly to $\big\{\mathcal{P}_{a,\alpha}\cdot \mathbf{1}_{[\mathcal{T}_x,\infty)}(s):s\in[0,1]\big\}$   in $\mathbf{D}([0,\infty),\mathbb{R})$ as $t\to\infty$.
 The convergence in the sense of finite-dimensional distributions follows directly from Theorem~\ref{Thm106}.
 Here we just need to prove the tightness.
 For any $0\leq r_1\leq r_2\leq r_3\leq 1$, we see that for $k=1,2$,
 \begin{equation*}
 Y^t_{r_{k+1}}-Y^t_{r_k} =\Big\{\begin{array}{ll}
 	t^{-1}\Delta \xi_{\mathcal{J}^{at}} , & \mathcal{J}^{at}\in(r_k,r_{k+1}];\cr
 	0, & \mbox{otherwise.}
 \end{array}
 \end{equation*}
 Hence  $|Y^t_{r_{2}}-Y^t_{r_1} |\cdot |Y^t_{r_{3}}-Y^t_{r_2} | \equiv 0$ a.s. and $\mathbf{P}\{ |Y^t_{r_{2}}-Y^t_{r_1} |\wedge |Y^t_{r_{3}}-Y^t_{r_2} | \geq \lambda \big|\tau_0>\mathcal{J}^{at} ,\xi_0=x\}=0$ for any $\lambda,t>0$.
 From Theorem~13.5 in \cite{Billingsley1999}, the sequence $\{Y^t_s:s\geq 0 \,\big|\,\tau_0>\mathcal{J}^{at},\xi_0=x\}_{t>0}$ is tight in $\mathbf{D}([0,\infty),\mathbb{R})$.
 \qed

 \begin{remark}\label{Remark01}
 For any $x\geq 0$, let $\tau_{-x}^d:= \inf\{k>0: \xi_k< -x\}$  and $\mathcal{J}_d^{x}:= \inf\{k\geq 1: \xi_k- \xi_{k-1}>x\}$.
 From Theorem~3.2 in \cite{DenisovShneer2013} we have as $n\to\infty$,
 \begin{equation}\label{eqn3.21}
 \mathbf{P}_x\{ \tau_0^d> n\} \sim \mathbf{E}_x[\tau_0^d]\cdot \mathbf{P}\{\xi_1> an\}.
 \end{equation}
 Following the previous argument, we also can establish the discrete versions of Theorem~\ref{Thm106} and \ref{Thm108} for $\{\xi_k:k=0,1,\cdots\}$ under Condition~\ref{Con2.1}.
 Here we show the results without detailed proofs.
 \begin{enumerate}
 \item[(1)] For any $x\geq 0$,  $b\geq a$ and $N>0$,
  \begin{equation}\label{eqn3.22}
  \mathbf{P}_x\big\{ \xi_{\mathcal{J}_d^{an}}- \xi_{\mathcal{J}_d^{an}-1}>bn, \mathcal{J}_d^{an} \leq N\,\big|\,\tau^d_0>n\big\}\to (b/a)^{-\alpha}\cdot \frac{\mathbf{E}_x[\tau^d_0\wedge N]}{\mathbf{E}_x[\tau^d_0]}.
  \end{equation}
 		
 \item[(2)]	$\{n^{-1}\xi_{[s]},s\geq 0\,|\,\tau^d_0>n, \xi_0=x\}$ converges weakly to $\{\mathcal{P}_{a,\alpha}\cdot \mathbf{1}_{\{\mathcal {T}^d_x<s\}},s\geq 0\}$  in $\mathbf{D}([0,\infty),\mathbb{R})$ as $t\to\infty$, where $\mathbf{P}\{\mathcal {T}^d_x\leq N\}=\mathbf{E}_x[\tau^d_0\wedge  N]/\mathbf{E}_x[\tau_0^d]$ for any $N\geq0$.
 \end{enumerate}
 	
 \end{remark}

 \subsection{Asymptotic results for conditional Laplace transforms}

 In this section, conditioned to $\xi_t>0$ for large $t$ or stay positive we provide limit theorems for the reflected processes together with several asymptotic results for the Laplace transforms of $\xi_t$.
 In the sequel of this section, we always assume both Condition~\ref{Con2.1} and the following condition hold.
 \begin{condition}\label{C2.2}
 For any $\delta>0$, we have as $x\to\infty$,
  \begin{equation*}
  \mathbf{P}\{ \xi_1 \in(x,x+\delta] \}= \frac{\alpha}{x}\mathbf{P}\{ \xi_1 >x \}\cdot\big[\delta +o(1)\big].
  \end{equation*}

 \end{condition}

 Note that this condition is in fact not really restrictive. For example, it holds in the following cases: (1) $\xi$ is a stable process with negative drift; (2) $\xi$ is a compound Poisson process with negative drift and Pareto distributed jumps; (3)  $\nu(dx):= \rho_\nu(x)dx$ with $\rho_\nu(x)\sim x^{-\alpha-1}\ell_0(x)$ as $x\to\infty$, where $\ell_0(x)$ is slowly varying.   Usually,  Condition~\ref{C2.2} holds if for any $\delta>0$  as $x\to\infty$,
 \begin{equation*}
 \nu(x,x+\delta]= \frac{\alpha}{x}\cdot \bar\nu(x)\cdot\big[\delta +o(1)\big].
 \end{equation*}
 Roughly speaking, like the argument below (\ref{eqn3.01}) we have as $x\to\infty$,
 \begin{eqnarray*}
 	\mathbf{P}\{ \xi_1 \in(x,x+\delta] \}\ar=\ar \mathbf{P}\Big\{\xi^1_1+ \int_0^1\int_{1}^\infty u N(ds,du) \in( x,x+\delta]\Big\} \cr
 	\ar\sim\ar \mathbf{P}\Big\{ \int_0^1\int_{1}^\infty u N(ds,du) \in( x,x+\delta]\Big\}\cr
 	\ar\ar\cr
 	\ar\sim\ar 1-\exp\big\{ \nu(x,x+\delta] \big\} \sim \nu(x,x+\delta].
 \end{eqnarray*}

 As a continuous analogue of Corollary~2.1 in \cite{DenisovDieker2008},
 the following lemma shows the asymptotic behavior of local probabilities for $\xi$ with proof given in Appendix.
 \begin{lemma}\label{Thm109}
 For any $\epsilon,\delta>0$, we have as $t\to\infty$,
 \begin{equation}\label{eqn3.24}
  \sup_{x\geq (\epsilon-a)t}\left|\frac{\mathbf{P}\{\xi_t\in[x,x+\delta)\}}{t\cdot \mathbf{P}\{\xi_1\in[at+x,at+x+\delta)\}}-1\right|\to 0.
 \end{equation}
 \end{lemma}

 \begin{proposition}\label{Thm110}
 For any $\lambda>0$ we have as $t\to\infty$,
 \begin{equation}\label{eqn3.25}
 \mathbf{E}[e^{-\lambda\xi_t};\xi_t\geq 0]\sim  \frac{\alpha}{a\lambda} \cdot\mathbf{P}\{\xi_1> at\}\quad\mbox{and}\quad
 \mathbf{E}[e^{\lambda\xi_t};\xi_t\leq 0]\sim  \frac{\alpha}{a\lambda} \cdot\mathbf{P}\{\xi_1> at\}.
 \end{equation}
 \end{proposition}
 \proof Here we just prove the first result and the second one can be proved similarly.
 For large $t$ we have
 \begin{equation*}
 \mathbf{E}[e^{-\lambda\xi_t};\xi_t\geq 0]= \mathbf{E}[e^{-\lambda\xi_t};\xi_t\in[ 0,\sqrt{t})]+o(e^{-\lambda \sqrt{t}}).
 \end{equation*}
 From Lemma~\ref{Thm109} and Condition~\ref{C2.2}, for any $n>1$ we have as $t\to\infty$,
  \begin{eqnarray*}
 \limsup_{t\to\infty} \frac{	\mathbf{E}[e^{-\lambda\xi_t};\xi_t\in[0,\sqrt{t})]}{\mathbf{P}\{\xi_1> at\}}
 \ar\leq\ar \limsup_{t\to\infty}\sum_{k=0}^{[n\sqrt{t}]}\int_{k/n}^{(k+1)/n} e^{-\lambda x}
  \frac{\mathbf{P}\{\xi_t\in dx\}}{\mathbf{P}\{\xi_1> at\}}\cr
  \ar\leq\ar  \sum_{k=0}^{\infty}e^{-\lambda \cdot k/n}\cdot
  \limsup_{t\to\infty}\frac{\mathbf{P}\{\xi_t\in [k/n,(k+1)/n)\}}{t\mathbf{P}\{ \xi_1\in[at+k/n,at+(k+1)/n) \}} \cr
  \ar\ar \quad \times   \limsup_{t\to\infty} \frac{t\mathbf{P}\{ \xi_1\in[at+k/n,at+(k+1)/n) \}}{\mathbf{P}\{\xi_1> at\}} \cr
  \ar\leq\ar \frac{\alpha}{a}\sum_{k=0}^{\infty}e^{-\lambda \cdot k/n} \cdot \frac{1}{n} ,
  	 \end{eqnarray*}
 which goes to $\alpha/(a\lambda)$ as $n\to\infty$. Similarly, we also can prove that as $n\to\infty$,
  \begin{eqnarray*}
 	\liminf_{t\to\infty} \frac{	\mathbf{E}[e^{-\lambda\xi_t};\xi_t\in[0,\sqrt{t})]}{\mathbf{P}\{\xi_1> at\}}
 	\ar\geq\ar \frac{\alpha}{a}\sum_{k=0}^{\infty}e^{-\lambda \cdot (k+1)/n} \cdot \frac{1}{n} \to \frac{\alpha}{a\lambda}.
 \end{eqnarray*}
 Putting these two estimates  together, we have as $t\to\infty$,
 \begin{eqnarray*}
 \mathbf{E}[e^{-\lambda\xi_t};\xi_t\geq 0]\sim  \mathbf{E}[e^{-\lambda\xi_t};\xi_t\in[0,\sqrt{t})]
 \sim \frac{\alpha}{a\lambda} \mathbf{P}\{\xi_1> at\}  .
 \end{eqnarray*}
 \qed

 As a preparation to study the asymptotic behavior of reflected processes, we provide the following useful proposition, which follows from Theorem~4(iii) in \cite{Cline1986}.
 \begin{proposition}\label{Thm111}
 Assume that $f(t)>0$ is regularly varying at $\infty$. For any two integrable functions $f_1,f_2 $ satisfying that $f_1(t)\sim c_1f(t)$ and $f_2(t)\sim c_2f(t)$ as $t\to\infty$ with $c_1,c_2\geq 0$, we have as $t\to\infty$,
 \begin{equation}\label{eqn3.26}
 \int_0^t f_1(t-s)f_2(s)ds\sim \Big( c_1\int_0^\infty f_2(s)ds+c_2\int_0^\infty f_1(s)ds \Big)\cdot f(t).
 \end{equation}
 \end{proposition}

 For any $z>0$ and $u,v> 0$, from Theorem~45.2 and 45.7 in  \cite{Sato1999} we have
 \begin{eqnarray}\label{eqn3.27}
 \lefteqn{\int_0^\infty ze^{-zt}\mathbf{E}\big[e^{-uS_t-v(S_t-\xi_t)}\big]dt}\ar\ar\cr
 \ar\ar=\exp\Big\{\int_0^\infty \frac{e^{-zt}}{t}\big(\mathbf{E}[e^{-u\xi_t};\xi_t\geq 0]+\mathbf{E}[e^{v\xi_t};\xi_t< 0]-1\big)dt\Big\}.
 \end{eqnarray}
 From the representations of $\kappa(z,u)$ and $\hat\kappa(z,u)$, we also have
 \begin{equation}\label{eqn3.28}
 \int_0^\infty ze^{-zt}\mathbf{E}\big[e^{-uS_t-v(S_t-\xi_t)}\big]dt
 =\frac{\kappa(z,0)}{\kappa(z,u)}\frac{\hat\kappa(z,0)}{\hat\kappa(z,v)}.
 \end{equation}
 Moreover, from Frullani's identity, we also have $\kappa(z,0)\hat\kappa(z,0) \sim  \mathcal{C}_0z$ as $z\to 0+$ with
 \begin{equation}\label{eqn3.29}
 \mathcal{C}_0:=\exp\Big\{-\int_0^\infty(1-e^{-t})\mathbf{P}\{\xi_t=0\}\frac{dt}{t}
 \Big\}.
 \end{equation}
 By the dominated convergence theorem,
 \begin{equation}\label{eqn3.30}
 \int_0^\infty \mathbf{E}[e^{-uS_t-v(S_t-\xi_t)}]dt
 = \lim_{z\to 0+}\frac{1}{z}\frac{\kappa(z,0)}{\kappa(z,u)}\frac{\hat\kappa(z,0)}{\hat\kappa(z,v)}
 =\frac{\mathcal{C}_0}{\kappa(0,u)\hat\kappa(0,v)}.
 \end{equation}
 From the identities $\int_0^\infty  e^{-ux}V(dx)= 1/\kappa(0,u) $ and $\int_0^\infty  e^{-vx}\hat{V}(dx)= 1/\hat\kappa(0,v) $; see (6) in \cite{Bertoin1996}[p.172] and then integration by parts, we have
 \begin{equation}\label{eqn3.31}
 \int_0^\infty \mathbf{E}\big[e^{-uS_t-v(S_t-\xi_t)}\big]dt
 =\mathcal{C}_0\int_0^\infty  e^{-ux}V(dx)\cdot \int_0^\infty  e^{-vy}\hat{V}(dy).
 \end{equation}
 We also can prove  the following result for $(I,I-\xi)$ in the same way,
 \begin{equation}\label{eqn3.32}
 \int_0^\infty\mathbf{E}\big[e^{uI_s+v(I_s-\xi_s)}\big]ds
 = \mathcal{C}_0\int_0^\infty e^{-ux}\hat{V}(dx)\cdot \int_0^\infty  e^{-vy}V(dy).
 \end{equation}
 Moreover, since Laplace transform is one-to-one, we  have for any $x,y\geq 0$,
 \begin{eqnarray}\label{eqn3.33}
 \int_0^\infty\mathbf{P}\{S_s\leq x, S_s-\xi_s\leq y\}ds\ar=\ar \mathcal{C}_0 V(x)\hat{V}(y),\label{eqn3.33.1}\\
 \int_0^\infty\mathbf{P}\{-I_s\leq x, \xi_s-I_s\leq y\}ds\ar=\ar \mathcal{C}_0 \hat{V}(x)V(y). \label{eqn3.33.2}
 \end{eqnarray}
 In the following two lemmas, we provide asymptotic results for the joint laws of the running supremum and infimum processes and their reflected processes.

 \begin{lemma}\label{Thm113}
 For any $u,v>0$, we have as $t\to\infty$,
 \begin{eqnarray}
 \frac{at}{\mathbf{P}\{\xi_1>at\}}\mathbf{E}\big[e^{-uS_t-v(S_t-\xi_t)}\big]\ar\to\ar \alpha\Big(\frac{1}{u}+\frac{1}{v}\Big)\int_0^\infty\mathbf{E}\big[e^{-uS_s-v(S_s-\xi_s)}\big]ds,\label{eqn3.34}\\
 \frac{at}{\mathbf{P}\{\xi_1> at\}}\mathbf{E}\big[e^{uI_t+v(I_t-\xi_t)}\big]
 \ar\to\ar \alpha\Big(\frac{1}{u}+\frac{1}{v}\Big)\int_0^\infty\mathbf{E}\big[e^{uI_s+v(I_s-\xi_s)}\big]ds.\label{eqn3.35}
 \end{eqnarray}
 \end{lemma}
 \proof Here we just prove the first result and the second one can be proved similarly. Taking the log on the both sides of (\ref{eqn3.27}) and then differentiating them with respect to $z$, we have
 \begin{eqnarray*}
 \lefteqn{\int_0^\infty (1-zt)e^{-zt}\mathbf{E}\big[e^{-uS_t-v(S_t-\xi_t)}\big]dt}\ar\ar\cr
 \ar=\ar-\int_0^\infty ze^{-zt}\mathbf{E}\big[e^{-uS_t-v(S_t-\xi_t)}\big]dt \int_0^\infty e^{-zs}\big(\mathbf{E}\big[e^{-u\xi_s};\xi_s\geq 0\big]+\mathbf{E}\big[e^{v\xi_s};\xi_s< 0\big]-1\big)ds,
 \end{eqnarray*}
which induces that
 \begin{eqnarray*}
	\lefteqn{\int_0^\infty  e^{-zt}\mathbf{E}\big[e^{-uS_t-v(S_t-\xi_t)}\big]dt-z\cdot\int_0^\infty e^{-zt}\cdot t\cdot \mathbf{E}\big[e^{-uS_t-v(S_t-\xi_t)}\big]dt}\ar\ar\cr
	\ar=\ar \int_0^\infty e^{-zt}\mathbf{E}[e^{-uS_t-v(S_t-\xi_t)}]dt \cr
	\ar\ar -z\cdot \int_0^\infty e^{-zt}\mathbf{E}[e^{-uS_t-v(S_t-\xi_t)}]dt\cdot \int_0^\infty e^{-zs}\big(\mathbf{E}[e^{-u\xi_s};\xi_s\geq 0]+\mathbf{E}[e^{v\xi_s};\xi_s< 0]\big)ds
\end{eqnarray*}
 and
 \begin{eqnarray*}
 \lefteqn{\int_0^\infty e^{-zt}\cdot t\cdot\mathbf{E}\big[e^{-uS_t-v(S_t-\xi_t)}\big]dt}\ar\ar\cr
 \ar=\ar\int_0^\infty e^{-zt}\mathbf{E}\big[e^{-uS_t-v(S_t-\xi_t)}\big]dt \cdot\int_0^\infty e^{-zs}\big(\mathbf{E}\big[e^{-u\xi_s};\xi_s\geq 0\big]+\mathbf{E}\big[e^{v\xi_s};\xi_s< 0\big]\big)ds\cr
 \ar=\ar \int_0^\infty e^{-zt}dt\int_0^t\mathbf{E}\big[e^{-uS_s-v(S_s-\xi_s)}\big]\big(\mathbf{E}\big[e^{-u\xi_{t-s}};\xi_{t-s}\geq 0\big]+\mathbf{E}\big[e^{v\xi_{t-s}};\xi_{t-s}< 0\big]\big)ds.
 \end{eqnarray*}
 Since Laplace transform is one-to-one, we have
 \begin{eqnarray}\label{eqn3.36}
 \lefteqn{t\cdot\mathbf{E}\big[e^{-uS_t-v(S_t-\xi_t)}\big]}\ar\ar \nonumber\\ \vspace{5pt}
 \ar=\ar\int_0^t\mathbf{E}\big[e^{-uS_s-v(S_s-\xi_s)}\big]\big(\mathbf{E}\big[e^{-u\xi_{t-s}};\xi_{t-s}\geq 0\big]+\mathbf{E}\big[e^{v\xi_{t-s}};\xi_{t-s}< 0\big]\big)ds.
 \end{eqnarray}
 Here we have $e^{-uS_t-v(S_t-\xi_t)}\leq e^{-uS_t}\leq e^{-u\xi_t} $ if $\xi_t\geq 0$ and $e^{-uS_t-v(S_t-\xi_t)}\leq e^{v\xi_t} $  if $\xi_t< 0$, which induce that
 \begin{eqnarray}\label{eqn3.34.1}
 	\lefteqn{t\cdot\mathbf{E}\big[e^{-uS_t-v(S_t-\xi_t)}\big]}
 	\ar\ar\cr
 	\ar\leq\ar
 	\int_0^t\mathbf{E}\big[e^{-u\xi_s};\xi_s\geq0\big]\Big(\mathbf{E}\big[e^{-u\xi_{t-s}};\xi_{t-s}\geq 0\big]+\mathbf{E}\big[e^{v\xi_{t-s}};\xi_{t-s}< 0\big]\Big)ds\cr
 	 \ar\ar + \int_0^t\mathbf{E}\big[e^{v\xi_s};\xi_s<0\big]\Big(\mathbf{E}\big[e^{-u\xi_{t-s}};\xi_{t-s}\geq 0\big]+\mathbf{E}\big[e^{v\xi_{t-s}};\xi_{t-s}< 0\big]\Big)ds.
 	\end{eqnarray}
 Moreover, from Proposition~\ref{Thm110} and~\ref{Thm111} we have as $t\to\infty$,
  \begin{eqnarray*}
  	\lefteqn{	\int_0^t\mathbf{E}\big[e^{-u\xi_s};\xi_s\geq0\big]\cdot \mathbf{E}\big[e^{-u\xi_{t-s}};\xi_{t-s}\geq 0\big] ds }\ar\ar\cr
  	\ar\sim\ar \frac{2\alpha}{au} \cdot	\int_0^\infty\mathbf{E}\big[e^{-u\xi_s};\xi_s\geq0\big] ds \cdot \mathbf{P}\{\xi_1> at\} \sim C \cdot \mathbf{P}\{\xi_1> at\}.
  	 \end{eqnarray*}
 Similar results for other terms in the right side of (\ref{eqn3.34.1}) can be proved in the same way. Putting these estimates together, there exists a constant $C>0$ such that for large $t$,
  \begin{eqnarray*}
  \mathbf{E}\big[e^{-uS_t-v(S_t-\xi_t)}\big]
 	\leq  \frac{C}{t} \mathbf{P}\{\xi_1> at\} = o\big(\mathbf{P}\{\xi_1> at\} \big).
  \end{eqnarray*}
 Applying Proposition~\ref{Thm110} and~\ref{Thm111} again to (\ref{eqn3.36}), we have as $t\to\infty$,
 \begin{equation*}
 t\cdot\mathbf{E}\big[e^{-uS_t-v(S_t-\xi_t)}\big]
 \sim \frac{\alpha}{a}\Big(\frac{1}{u}+\frac{1}{v}\Big)\int_0^\infty
 \mathbf{E}\big[e^{-uS_s-v(S_s-\xi_s)}\big]ds \cdot \mathbf{P}\{\xi_1> at\}.
 \end{equation*}
 Here we have gotten the desired result.
 \qed

 \begin{lemma}\label{Thm114}
 For any $x,y\geq 0$, we have as $t\to\infty$,
 \begin{eqnarray*}
 \frac{at}{\mathbf{P}\{\xi_1> at\}}\mathbf{P}\big\{S_t\leq x,S_t-\xi_t\leq y\big\}\ar\to\ar \mathcal{C}_0\alpha\Big[\hat{V}(y)
 \int_0^xV(z)dz+V(x)\int_0^y\hat{V}(z)dz\Big],\label{eqn3.37}\cr
 \frac{at}{\mathbf{P}\{\xi_1> at\}}\mathbf{P}\big\{-I_t\leq x,\xi_t- I_t\leq y\big\}\ar\to\ar \mathcal{C}_0\alpha\Big[ V(y)
 \int_0^x\hat{V}(z)dz+\hat{V}(x)\int_0^yV(z)dz\Big].
 \end{eqnarray*}
 \end{lemma}
 \proof  From (\ref{eqn3.36}), we have
 \begin{eqnarray*}
  \lefteqn{t\cdot\int_0^\infty \int_0^\infty e^{-ux-vy} \mathbf{P}\{S_t\in dx, S_t-\xi_t\in dy \} }\ar\ar \nonumber\\ \vspace{5pt}
  \ar=\ar\int_0^t ds\int_0^\infty \int_0^\infty  \int_0^\infty e^{-u(x_1+z)-vy_1} \mathbf{P}\{S_s\in dx_1, S_s-\xi_s\in dy_1 \} \mathbf{P}\{ \xi_{t-s}\in dz \} \cr
  \ar\ar +\int_0^t ds\int_0^\infty \int_0^\infty \int_0^\infty e^{-ux_1-v(y_1+z)} \mathbf{P}\{S_s\in dx_1, S_s-\xi_s\in dy_1 \}    \mathbf{P}\{ -\xi_{t-s}\in dz \}  .
  \end{eqnarray*}
 From this and the one-to-one property of Laplace transform,
 \begin{eqnarray}\label{eqn3.39}
 t\cdot\mathbf{P}\{S_t\leq x,S_t-\xi_t\leq y\}
  \ar=\ar\int_0^tds \int_0^x \mathbf{P}\{S_s\leq x-z,S_s-\xi_s\leq y\}\mathbf{P}\{\xi_{t-s}\in dz\}\cr
 \ar\ar +\int_0^tds \int_0^y \mathbf{P}\{S_s\leq x,S_s-\xi_s\leq y-z\}\mathbf{P}\{-\xi_{t-s}\in dz\}.\quad\quad
 \end{eqnarray}
 It is easy to see that
  \begin{eqnarray*}
  	t\cdot \mathbf{P}\{S_t\leq x,S_t-\xi_t\leq y\}\leq \int_0^t\mathbf{P}\{S_s\leq x,S_s-\xi_s\leq y\}\mathbf{P}\{-y\leq\xi_{t-s}\leq x\}ds,
  	 \end{eqnarray*}
 where the first probability on the right side of above inequality can be bounded as follows
   \begin{eqnarray*}
   	\mathbf{P}\{S_s\leq x,S_s-\xi_s\leq y\}
   	 \leq \mathbf{P}\{0\leq \xi_s\leq x \} +\mathbf{P}\{-y\leq \xi_s\leq 0\}
   	 = \mathbf{P}\{-y\leq \xi_s\leq x\} .
   \end{eqnarray*}
 From  Lemma~\ref{Thm109} and Condition~\ref{C2.2}, we have
 \begin{equation}\label{eqn3.40}
 \mathbf{P}\{-y\leq\xi_{t}\leq x\}\sim \frac{\alpha (x+y)}{a}\mathbf{P}\{\xi_1> at\}.
 \end{equation}
 Putting these estimates together, from  Proposition~\ref{Thm111} there exists a constant $C>0$ such that for large $t$,
 \begin{eqnarray*}
 	\mathbf{P}\{S_t\leq x,S_t-\xi_t\leq y\}\ar\leq\ar\frac{C}{t} \mathbf{P}\{\xi_1>at\} = o\big( \mathbf{P}\{\xi_1>at\}  \big).
 \end{eqnarray*}
 Applying Proposition~\ref{Thm111} with this result and (\ref{eqn3.40}), for any $n>1$ we have as $t\to\infty$,
\begin{eqnarray*}
	\lefteqn{\frac{a}{\mathbf{P}\{\xi_1>at\}}\int_0^{t}ds\int_0^x \mathbf{P}\{S_s\leq x-z,S_s-\xi_s\leq y\}\mathbf{P}\{\xi_{t-s}\in dz\}} \ar\ar\cr
	\ar\leq\ar \sum_{k=0}^{[nx]}\frac{a}{\mathbf{P}\{\xi_1>at\}}\int_0^{ t} \mathbf{P}\{S_s\leq x-k/n,S_s-\xi_s\leq y\} \cdot \mathbf{P}\{\xi_{t-s}\in [k/n,(k+1)/n)\}ds\cr
	\ar\to \ar \sum_{k=0}^{[nx]}\frac{\alpha}{n}\int_0^\infty \mathbf{P}\{S_s\leq x-k/n,S_s-\xi_s\leq y\} ds.
	\end{eqnarray*}
 Let $n\to\infty$,  we have
 \begin{eqnarray*}
\lefteqn{ \limsup_{t\to\infty} \frac{a}{\mathbf{P}\{\xi_1>at\}}\int_0^tds\int_0^x \mathbf{P}\{S_s\leq x-z,S_s-\xi_s\leq y\}\mathbf{P}\{\xi_{t-s}\in dz\}} \ar\ar\cr
\ar\ar\leq  \alpha \int_0^\infty ds\int_0^x\mathbf{P}\{S_s\leq x-z,S_s-\xi_s\leq y\}dz
= \mathcal{C}_0\alpha\hat{V}(y) \int_0^xV(z)dz.
 \end{eqnarray*}
 Similarly, we also have
 \begin{eqnarray*}
 \lefteqn{ \liminf_{t\to\infty} \frac{a}{\mathbf{P}\{\xi_1>at\}}\int_0^tds\int_0^x \mathbf{P}\{S_s\leq x-z,S_s-\xi_s\leq y\}\mathbf{P}\{\xi_{t-s}\in dz\}} \ar\ar\cr
 \ar\ar\geq  \alpha \int_0^\infty ds\int_0^x\mathbf{P}\{S_s\leq x-z,S_s-\xi_s\leq y\}dz
 = \mathcal{C}_0\alpha\hat{V}(y) \int_0^xV(z)dz.
 \end{eqnarray*}
 Putting these two estimate together, we have
 \begin{equation*}
 \frac{a}{\mathbf{P}\{\xi_1> at\}}\int_0^{t}ds\int_0^x \mathbf{P}\{S_s\leq x-z,S_s-\xi_s\leq y\}\mathbf{P}\{\xi_{t-s}\in dz\}\to \mathcal{C}_0\alpha\hat{V}(y)
 \int_0^xV(z)dz.
 \end{equation*}
 Similarly, we also have
 \begin{equation*}
 \frac{a}{\mathbf{P}\{\xi_1>at\}}\int_0^{t}ds\int_0^y \mathbf{P}\{S_s\leq x,S_s-\xi_s\leq y-z\}\mathbf{P}\{-\xi_{t-s}\in dz\}\to\mathcal{C}_0\alpha V(x)
 \int_0^y\hat{V}(z)dz.
 \end{equation*}
 Taking these results back into (\ref{eqn3.39}), we can get the first desired result.
 The second one can be proved similarly.
 \qed

 \begin{theorem}\label{Thm115}
 For any $x, \lambda>0$, we have as $t\to\infty$,
 \begin{eqnarray}
 \frac{at}{\mathbf{P}\{\xi_1> at\}}\mathbf{E}_x[e^{-\lambda\xi_t};\tau_0>t]
 \ar\to\ar \mathcal{C}_0\alpha \cdot    \hat{V}(x)\int_0^\infty e^{-\lambda y}V(y)dy ,\label{eqn3.42}\\
  \frac{at}{\mathbf{P}\{\xi_1>at\}}\mathbf{E}_{-x}[e^{\lambda\xi_t};\tau_0^+>t]
 \ar \to\ar  \mathcal{C}_0\alpha \cdot  V(x)\int_0^\infty e^{-\lambda y}\hat{V}(y)dy . \label{eqn3.43}
 \end{eqnarray}
 \end{theorem}
 \proof Here we still just prove the first statement and the second one can be proved similarly. Applying duality of $\xi$ to the second equality below, we have
 \begin{eqnarray*}
\mathbf{E}_x\big[e^{-\lambda\xi_t};\tau_0>t\big]
 =  e^{\lambda x}\mathbf{E}\big[e^{-\lambda\xi_t};I_t\geq -x\big]
 \ar=\ar   e^{-\lambda x}\mathbf{E}\Big[e^{-\lambda(\xi_t-\xi_0)};\Big(\xi_t-\inf_{s\leq t}\xi_s\Big)-\xi_t\leq x\Big]\cr
 \ar=\ar   e^{-\lambda x} \mathbf{E}\Big[e^{-\lambda[S_t-(S_t-\xi_t)]};S_t-\xi_t\leq x\Big]\cr
 \ar=\ar  e^{-\lambda x}\int_0^\infty\int_0^xe^{-\lambda(y-z)}\mathbf{P}\big\{S_t\in  dy, S_t-\xi_t\in dz\big\}.
 \end{eqnarray*}
 For $k>8x$ and $\theta\in(\lambda/2,\lambda)$, we have
 \begin{eqnarray*}
 \lefteqn{\int_k^\infty\int_0^xe^{-\lambda(y-z)}\mathbf{P}\big\{S_t\in dy, S_t-\xi_t\in dz\big\}}\ar\ar\cr
 \ar=\ar \int_k^\infty\int_0^xe^{-\frac{\lambda}{2} y-\lambda(y/2-z)}\mathbf{P}\big\{S_t\in dy, S_t-\xi_t\in dz\big\}\cr
 \ar\leq\ar  e^{-\frac{\theta}{2} k}\int_0^\infty\int_0^x e^{-\frac{\lambda-\theta}{2} y-(\frac{\theta}{2} y-\lambda z)}\mathbf{P}\big\{S_t\in dy, S_t-\xi_t\in dz\big\}\cr
 \ar\leq\ar  e^{-\frac{\theta}{2} k}\int_0^\infty\int_0^\infty e^{-\frac{\lambda-\theta}{2} y-\lambda z}\mathbf{P}\big\{S_t\in dy, S_t-\xi_t\in dz\big\}
 =e^{-\frac{\theta}{2} k}\mathbf{E}\big[e^{-\frac{\lambda-\theta}{2}S_t-\lambda ( S_t-\xi_t)}\big].
 \end{eqnarray*}
Here the second inequality follows from the fact that $\frac{\theta}{2} y-\lambda z\geq 4\theta x -\lambda z  \geq \lambda z$ for any $z\in[0,x]$ and $y\geq 8x$.
 From this and  Lemma~\ref{Thm113},
 \begin{equation}
 \lim_{k\to\infty}\limsup_{t\to\infty}\frac{at}{\mathbf{P}\{\xi_1> at\}}\int_k^\infty\int_0^xe^{-\lambda(y-z)}\mathbf{P}\big\{S_t\in dy, S_t-\xi_t\in dz\big\}=0.
 \end{equation}
 From this estimate and Lemma~\ref{Thm114}, we have
 \begin{eqnarray*}
 \lefteqn{\lim_{t\to\infty}\frac{at}{\mathbf{P}\{\xi_1> at\}}\mathbf{E}[e^{-\lambda\xi_t};I_t\geq -x]}\ar\ar\cr
 \ar=\ar \lim_{k\to\infty }\lim_{t\to\infty} \frac{at}{\mathbf{P}\{\xi_1> at\}} \int_0^k\int_0^xe^{-\lambda(y-z)}\mathbf{P}\{S_t\in dy, S_t-\xi_t\in dz\}\cr
 \ar=\ar C_0\alpha\int_0^\infty \int_0^xe^{-\lambda(y-z)}\cdot d\Big[\hat{V}(z)
 \int_0^yV(w)dw+V(y)
 \int_0^z\hat{V}(w)dw\Big]\cr
 \ar=\ar C_0\alpha\int_0^\infty e^{-\lambda y}V(y)dy\int_0^xe^{\lambda z}d\hat{V}(z)
 +C_0\alpha\int_0^\infty e^{-\lambda y} dV(y)\int_0^xe^{\lambda z}
 \hat{V}(z)dz.
 \end{eqnarray*}
 By integration by parts, we can immediately get the desired result.
 \qed

 \begin{remark}
  Following the argument above, we also can prove the analogue of Theorem~\ref{Thm115} for the random walk $\{\xi_n:n=0,1,\cdots\}$ under Condition~\ref{Con2.1} and \ref{C2.2}, i.e. let $\tau_0^{d+}:= \inf\{n> 0: \xi_n> 0  \}$,  for any $x\geq 0$ and $\lambda>0$ there exist  constants $C_1,C_2>0$ such that as $n\to\infty$,
  \begin{equation}\label{eqn3.44}
  \mathbf{E}_x[e^{-\lambda\xi_n};\tau_0^d>n]
  \sim C_1 \frac{\mathbf{P}\{\xi_1> an\}}{an}
  \quad \mbox{and}\quad
  \mathbf{E}_x[e^{\lambda\xi_n};\tau_0^{d+}>n]
  \sim C_2 \frac{\mathbf{P}\{\xi_1> an\}}{an}.
  \end{equation}
 \end{remark}


 \section{Asymptotic results for exponential functionals}
 \setcounter{equation}{0}

  In this section, we study the asymptotic behaviors of exponential functionals (\ref{eq1.1}) of heavy-tailed L\'evy processes with the help of conditional limit results introduced in the last section.
  From Theorem~1 in \cite{BertoinYor2005}, $\mathbf{P}\{A_\infty(\xi)< \infty\} = 1$ if and only if $\mathbf{P}\{A_\infty(\xi)< \infty\}>0$; equivalently, if and only if $a< 0$.
 For any $x>0$ and $t\geq 0$, we have
 \begin{equation}\label{eqn4.01}
 (1-e^{-t})\mathbf{P}\{A_{t}(\xi)\leq x\} \leq   \int_0^\infty e^{-s}\mathbf{P}\{A_s(\xi)\leq x\}ds = \mathbf{P}\{A_{\bf e}(\xi)\leq x\} ,
 \end{equation}
 where ${\bf e}$ is an exponentially distributed random variable with parameter $1$ and independent of $\xi$.
 By Theorem~2.19 in \cite{PatieSavov2018}, we have $\mathbf{P}\{A_{\bf e}(\xi)\leq x\}\sim C x$ as $x\to 0+$ and hence there exists a constant $C_t>0$ such that for any $x\geq 0$,
 \begin{equation}\label{eqn4.02}
 \mathbf{P}\{A_t(\xi)\leq x\}\leq \frac{e^t}{e^t-1}\mathbf{P}\{A_{\bf e}(\xi)\leq x\}\leq  C_t x.
 \end{equation}
 Moreover, they also proved that $ \mathbf{E}\big[|A_t(\xi)|^{-\kappa}\big]<\infty$ for any  $\kappa\in(0,1)$ and $t>0$; see Theorem~2.18 in \cite{PatieSavov2018}.

 We now start to study the asymptotic behavior of the expectation $\mathbf{E}[F(A_t(\xi))]$ as $t\to\infty$ for some function $F$ defined on $(0,\infty)$.
 To simplify the presentation of the results, we also assume that $F$ satisfies the following two conditions.
 \begin{condition}\label{ConF1}
 	$F$ is bounded, positive, non-increasing and $K_0:=\sup_{x>0}x^{\beta}F(x)<\infty$ for some $\beta\in(0,1)$.
 \end{condition}
 \begin{condition}\label{ConF2}
 	For any $\delta>0$, there exists a constant $K_\delta>0$ such that $|F(x)-F(y)|\leq K_\delta |x-y|$ for any $x,y\geq\delta$.
 \end{condition}

 If  $a<0$, we have $\mathbf{P}\{A_\infty(\xi)< \infty\} = 1$ and hence $\mathbf{E}[F(A_t(\xi))]\to \mathbf{E}[F(A_\infty(\xi))]<\infty$ as $t\to\infty$.
 If $a=0$ and $\rho_0:=\lim_{t\to\infty}\mathbf{P}\{\xi_t<0\}\in[0,1)$, we have $ \mathbf{E}[F(A_t(\xi))]\sim  t^{-\rho_0}\ell_0(t)$ as $t\to\infty$, where  $\ell_0(t)$ is a slowly varying function at $\infty$; see Theorem~2.18(2) in \cite{PatieSavov2018}\footnote{Their result holds for any $F$ satisfying that  $\sup_{x>0}x^{\beta}F(x)<\infty$ for some $\beta\in(0,1)$.}
 In this section, we consider the asymptotic behavior of $\mathbf{E}[F(A_t(\xi))]$ with $a>0$ and $\xi$ satisfying Condition~\ref{Con2.1} and \ref{C2.2}. Again, let $\xi'$ be an independent copy of $\xi$.
 Let $\mathbf{J}$ be an $\mathbb{R}_+$-valued random variable independent of $\xi$ and $\xi'$ with $\mathbf{P}\{ \mathbf{J} >x  \}= \bar\nu(x)$ for large $x$.
 Actually, all the following results hold for any $\mathbf{J}$ satisfying that $\mathbf{P}\{ \mathbf{J} >x  \}\sim \bar\nu(x)$ as $x\to\infty$.  For any function $g\in \mathbf{D}([0,\infty),\mathbb{R})$ and $s\geq 0$, define
  \begin{equation*}
 C_F(g,s):= \lim_{t\to\infty}\mathbf{E} \big[F\big(A_s(g)+e^{-g_{s}-\mathbf{J}}A_{t-s}(\xi')\big)\big|\mathbf{J}>at\big],
 \end{equation*}
 where $A_s(g)$ is defined as (\ref{eq1.1}) with $\xi$ replaced by $g$. The existence of the limit above will be proved in Lemma~\ref{Thm408}.
 From the independence between $\xi$ and $\xi'$, we see that the stochastic processes $\{C_F(\xi,s):s\geq 0\}$ is well defined.
 The main theorem of this section is the following:

 \begin{theorem}\label{T1}
 Assume that $a>0$ and Condition~\ref{Con2.1} and \ref{C2.2} hold, we have the finite and nonzero limit
 \begin{equation}\label{eqn4.03}
 \lim_{t\to\infty}\frac{\mathbf{E}[F(A_t(\xi))]}{\mathbf{P}\{\xi_1>at\}} =\int_0^\infty \mathbf{E}\big[ C_F(\xi,s)\big]ds<\infty,
 \end{equation}
 \end{theorem}

 Before showing the proof for this theorem, in the following lemma we study the effect of the initial state $\xi_0=-x$ on the expectation $\mathbf{E}_{-x}[F(A_t(\xi))]=\mathbf{E}\big[F(e^{x }A_t(\xi))\big]$.
 This offers us a criticality to identify the sample paths that make main contribution to the expectation (\ref{eq1.2}).

 \begin{lemma}\label{Thm403}
 For any $t\geq 0$ and $q>0$, there exist two constants $C,\lambda_0>0$ such that for any $x\in\mathbb{R}$,
 \begin{equation}\label{eqn4.04}
 \mathbf{E}\big[\big|F(e^{x }A_t(\xi))\big|^q\big]+\mathbf{E}\big[F(e^{ x}A_t(\xi))\cdot \xi_t^+\big]\leq C e^{-\lambda_0 x}.
 \end{equation}
 \end{lemma}
 \proof Here we just prove this result with $t=1$ and $x\geq 0$. Other cases can be proved in the same way.
 Firstly,
 \begin{eqnarray*}
 \lefteqn{\mathbf{E}[|F(e^{x}A_1(\xi))|^q]
 =\int_{-\infty}^\infty  |F(e^{ x-y})|^q \, d\mathbf{P}\{-\log A_1(\xi)\leq y\}}\ar\ar\cr
 \ar=\ar   \int_{-\infty}^{x/2} |F(e^{ x-y})|^q \, d\mathbf{P}\{A_1(\xi)\geq e^{-y}\} + \int_{x/2}^\infty  |F(e^{ x-y})|^q \, d\mathbf{P}\{A_1(\xi)\geq e^{-y}\} .
 \end{eqnarray*}
 Since $F $ is non-increasing, we have
 \begin{equation*}
 \int_{-\infty}^{x/2} |F(e^{ x-y})|^q \,  d\mathbf{P}\{A_1(\xi)\geq e^{-y}\} \leq |F(e^{ x/2})|^q \leq Ce^{-q\beta x/2}.
 \end{equation*}
 From the boundedness of $F$ and (\ref{eqn4.02}), we also have
 \begin{equation*}
 \int_{x/2}^\infty  |F(e^{ x-y})|^qd\mathbf{P}\{A_1(\xi)\geq e^{-y}\}
 \leq  C\mathbf{P}\{A_1(\xi)\leq  e^{- x/2}\}<C e^{- x/2}.
 \end{equation*}
 Putting all estimates above together, we have
 $\mathbf{E}[|F(e^{x }A_1(\xi))|^q]\leq Ce^{-(q\beta \wedge 1)x/2}$.
 By H\"older's inequality, for any $p,q> 1$ satisfying that $p< \alpha$ and $\frac{1}{p}+\frac{1}{q}=1$ we have
 \begin{equation*}
 \mathbf{E}\big[F(e^{ x}A_1(\xi))\cdot \xi_1^+\big]\leq \mathbf{E}\big[\big|F(e^{ x}A_1(\xi))\big|^q\big]^{1/q}\cdot\mathbf{E}\big[ |\xi_1^+|^p\big]^{1/p}\leq C e^{-\lambda_0x}
 \end{equation*}
 with $\lambda_0:= (q\beta \wedge 1)/(2q)$.
 \qed

 This lemma shows that the asymptotic of the expectation (\ref{eq1.2}) may mainly depend on the sample paths with slowly decreasing local infimum.
 From Lemma~\ref{Thm107} and Theorem~\ref{Thm108}, we see that the sample paths decrease slowly if and only if there is a large jump that occurs early.
 To show clearly the main ideas of the proof for Theorem~\ref{T1}, we write it into several steps with a series of lemmas.
 In Section~4.1, we prove that the contribution of sample paths with late arrival of large jump to the expectation (\ref{eq1.2}) is insignificant,  i.e.,  for $t,N>0$ large enough,
 \begin{equation}\label{eqn4.05}
 \mathbf{E}[F(A_t(\xi));\mathcal{J}^{at}>N]=o(\mathbf{P}\{\xi_1>at\}).
 \end{equation}
 In Section~4.2, we analyze the exact contribution of sample paths with early arrival of large jump to the expectation.
 Enlightened by the conditional asymptotic behaviors of $\xi$ provided in  Theorem~\ref{Thm108},  we observe that $A_t(\xi)$ increases very slowly after the arrival of large jump, which results that $F(A_t(\xi))$ decreases so slowly that it can be well approximated by $F(A_T(\xi))$ for large $T$,  i.e., for large $N$,
 \begin{equation*}
 \mathbf{E}\big[F(A_t(\xi));\mathcal{J}^{at}<N]
 \sim \mathbf{E}\big[F( A_{T}(\xi));\mathcal{J}^{at}<N\big]\sim C(T) \mathbf{P}\{\xi_1>at\}.
 \end{equation*}
 Based on these estimates, in Section~4.3 we give the proof for Theorem~\ref{T1}.

 \subsection{Contribution of sample paths with late arrival of large jump}

 In this section, we prove that the contribution of sample paths with late arrival of large jump to the expectation (\ref{eq1.2}) can be ignored.
 Recall  $\mathcal{J}^{x}_d$ defined in Remark~\ref{Remark01}.
 From the fact that $\mathbf{P}\{\mathcal{N}_1^{ax}>0\} \sim \mathbf{P}\{\xi_1\geq ax\}$ as $x\to\infty$; see (\ref{eqn3.01}) and (\ref{eqn3.04}), for large $x$ the following proposition shows that $\mathcal{J}^{ax}>N$ if and only if $\mathcal{J}^{ax}_d>N$.

 \begin{proposition}\label{Thm404}
 For any $N\geq 1$,	we have as $x\to\infty$,
 \begin{equation}\label{eqn4.07}
 \frac{\mathbf{P}\big\{\{\mathcal{N}_1^{ax}>0\} \Delta \{  \xi_1>ax \}  \big\} }{\mathbf{P}\{\xi_1>ax\} }\to 0
 \quad\mbox{and}\quad
 \frac{\mathbf{P}\big\{\{[\mathcal{J}^{ax}]\geq N\} \Delta \{ \mathcal{J}^{ax}_d\geq N\} \big\}}{\mathbf{P}\{\xi_1> ax\}} \to 0.
 \end{equation}
 \end{proposition}
 \proof  From  the fact that $\mathbf{P}\{\mathcal{N}_1^{ax}>0\} \sim \mathbf{P}\{\xi_1\geq ax\}$ as $x\to\infty$, we have
 \begin{eqnarray*}
 	\frac{   \mathbf{P}\big\{\{\mathcal{N}_1^{ax}>0\} \Delta \{  \xi_1> ax \}  \big\} }{\mathbf{P}\{\xi_1> ax\} }\ar=\ar 1- \frac{   \mathbf{P}\{\mathcal{N}_1^{ax}>0\}   }{\mathbf{P}\{\xi_1> ax\} }+ \frac{ 2 \mathbf{P}\{\xi_1\leq ax, \mathcal{N}^{ax}_1>0  \}}{\mathbf{P}\{\xi_1> ax\}}\cr
 	\ar\sim\ar \frac{ 2 \mathbf{P}\{\xi_1\leq ax, \mathcal{N}^{ax}_1>0  \}}{\mathbf{P}\{\xi_1> ax\}}.
 \end{eqnarray*}
 For any $b>a$,
 \begin{equation*}
 \mathbf{P}\big\{\xi_1\leq ax, \mathcal{N}^{ax}_1>0\big\}
 \leq \mathbf{P}\big\{\xi_1\leq ax,\mathcal{N}^{bx}_1>0\big\}+ \mathbf{P}\big\{\mathcal{N}^{ax}_1-\mathcal{N}^{bx}_1>0\big\}.
 \end{equation*}
 From  (\ref{eqn3.03}) we have $\mathcal{N}^{ax}_1-\mathcal{N}^{bx}_1= \int_0^1 \int_{ax}^{bx} N(ds,du)$ and as $x\to\infty$,
  \begin{equation*}
 \mathbf{P}\big\{\mathcal{N}^{ax}_1-\mathcal{N}^{bx}_1>0\big\}= 1-\exp\big\{  -\big[\bar\nu(ax)-\bar\nu(bx)\big] \big\}
 \sim \bar\nu(ax)-\bar\nu(bx),
  \end{equation*}
 which immediately induces that
 \begin{equation*}
 \frac{ \mathbf{P}\{   \mathcal{N}^{ax}_1-\mathcal{N}^{bx}_1>0  \}}{\mathbf{P}\{\xi_1> ax\}} \to 1-(a/b)^{\alpha}.
 \end{equation*}
 Moreover, since $\xi_1^1$ is independent of $\mathcal{N}^{bx}_1$ for $bx>1$, then
 \begin{eqnarray*}
 \mathbf{P}\big\{  \xi_1\leq ax, \mathcal{N}^{bx}_1>0  \big\}
 \ar\leq\ar \mathbf{P} \big\{  \xi^1_1+bx\leq ax, \mathcal{N}^{bx}_1>0  \big\} \cr
 \ar\ar\cr
 \ar=\ar\mathbf{P}\big\{  \xi^1_1\leq (a-b)x \big\} \cdot \mathbf{P}\big\{ \mathcal{N}^{bx}_1>0  \big\} = o\big(\mathbf{P}\{ \mathcal{N}^{bx}_1>0  \} \big).
 \end{eqnarray*}
 Putting all results above together, we have
 \begin{equation*}
 \limsup_{x\to\infty}
 \frac{   \mathbf{P}\big\{\{\mathcal{N}_1^{ax}>0\} \Delta \{  \xi_1> ax \}  \big\} }{\mathbf{P}\{\xi_1> ax\} }\leq  1-(a/b)^{\alpha},
 \end{equation*}
 which vanishes as $b\to a+$.
 We start to prove the second result.
 By the independent increments of $\xi$, we have for any $k<N$,
 \begin{eqnarray*}
 \lefteqn{\mathbf{P}\Big\{ [\mathcal{J}^{ax}]\geq N, \mathcal{J}^{ax}_d=k \Big\}
  = \prod_{i=1}^{k-1} \mathbf{P}\Big\{ \int_{i-1}^{i}\int_{ax}^\infty N(ds,du)=0, \xi_{i}-\xi_{i-1}\leq ax \Big\}}\ar\ar\cr
 \ar\ar \times \mathbf{P}\Big\{ \int_{k-1}^{k}\int_{ax}^\infty N(ds,du)=0, \xi_{k}-\xi_{k-1}> ax \Big\}  \times  \mathbf{P}\Big\{ \int_{k}^{N}\int_{ax}^\infty N(ds,du)=0\Big\}\cr
 \ar\ar\cr
 \ar=\ar \big|\mathbf{P}\{ \mathcal{N}_1^{ax}=0, \xi_1\leq ax \}\big|^{k-1} \cdot \mathbf{P}\big\{ \mathcal{N}_1^{ax}=0, \xi_1> ax \big\} \cdot \big|\mathbf{P}\{ \mathcal{N}_1^{ax}=0\}\big|^{N-k}.
 \end{eqnarray*}
 Both the first and the last probability on the right side of the last equality go to $1$ as $x\to\infty$, which immediately induces that
 \begin{equation*}
 \mathbf{P}\big\{ [\mathcal{J}^{ax}]\geq N, \mathcal{J}^{ax}_d< N \big\}
 \sim  (N-1)\cdot\mathbf{P}\big\{ \mathcal{N}_1^{ax}=0, \xi_1> ax \big\}.
 \end{equation*}
 Similarly, we also have as $x\to\infty$,
 \begin{equation*}
  \mathbf{P}\big\{ [\mathcal{J}^{ax}]< N, \mathcal{J}^{ax}_d\geq  N \big\}
  \sim  (N-1)\cdot\mathbf{P}\big\{ \mathcal{N}_1^{ax}>0, \xi_1\leq ax \big\}  .
 \end{equation*}
 Putting these two estimates together, we have as $x\to\infty$,
 \begin{equation*}
 \mathbf{P}\{\{[\mathcal{J}^{ax}]\geq N\} \Delta \{ \mathcal{J}^{ax}_d\geq N\} \}
 \sim (N-1)\cdot \mathbf{P}\{\{\mathcal{N}_1^{ax}>0\} \Delta \{  \xi_1> ax \}  \}
 \end{equation*}
 and the second result follows from this and the first result in (\ref{eqn4.07}).
 \qed

 From the last proposition, it suffices to prove that $\mathbf{E}[F(A_t(\xi));\mathcal{J}_d^{at}>N]=o(\mathbf{P}\{\xi_1>at\})$ for large $t$ and $N$.
 For any $n\geq 1$, let $I^d_n :=\inf\{\xi_k:k=0,1\cdots n\}$ and $ g^d_n:= \sup\{0\leq k\leq n:\xi_k=I^d_n \}$.
 According to the distance between the large jump and the local infimum, we split the expectation $\mathbf{E}[F(A_t(\xi));\mathcal{J}_d^{at}>N]$ into the following two terms: for $K>0$,
 \begin{eqnarray}\label{eqn4.08}
  \mathbf{E}[F(A_t(\xi));\mathcal{J}_d^{at}>N, g^d_{[t]}<K]+ \mathbf{E}[F(A_t(\xi));\mathcal{J}_d^{at}>N;g^d_{[t]}\geq K].
 \end{eqnarray}
 The following lemma shows that the local infimum mostly is close to the first large jump, which immediately induces that the first term in (\ref{eqn4.08}) can be ignored.

 \begin{lemma}\label{Thm411}
 For any fixed $K>0$, we have
 \begin{equation}\label{eqn4.09}
 \lim_{N\rightarrow \infty}\limsup_{n\rightarrow \infty}\frac{\mathbf{P}\{\mathcal{J}^{an}>N,g^d_n<K\}}{\mathbf{P}\{\xi_1>an\}}=0.
 \end{equation}
 \end{lemma}
 \proof From Proposition~\ref{Thm404}, it suffices to prove the following result with $K\in\mathbb{Z}_+$:
 \begin{equation*}
 \lim_{N\rightarrow \infty}\limsup_{n\rightarrow \infty}\frac{\mathbf{P}\{\mathcal{J}_d^{an}>N,g^d_n<K\}}{\mathbf{P}\{\xi_1>an\}}=0,
 \end{equation*}
 By the duality of $\xi$,  for $n>N>K$ we have
 \begin{eqnarray*}
 \mathbf{P}\{\mathcal{J}_d^{an}>N,g^d_n<K\}
 \ar=\ar\sum_{k=0}^{K-1} \mathbf{P}\{\mathcal{J}_d^{an}>N,g^d_n=k\}\cr
 \ar\leq\ar \sum_{k=0}^{K-1} \mathbf{P}\{\mathcal{J}_d^{an}>N-k,\tau_0^d\geq n-k\} \mathbf{P}\{ g_k^d= k \}\cr
 \ar\leq \ar \mathbf{P}\{ \tau_0^d\geq n-K\} \cdot K\cdot \mathbf{P}\big\{\mathcal{J}_d^{an}>N-K|\tau_0^d\geq n-K\big\}.
 \end{eqnarray*}
 From  (\ref{eqn3.21}) and (\ref{eqn3.22}), we have
 \begin{equation*}
 \lim_{n\to\infty}\frac{ \mathbf{P}\big\{\mathcal{J}_d^{an}>N,g^d_n<K\big\}}{\mathbf{P}\{\xi_1>an\}}
 \leq C\cdot K\cdot  \frac{\mathbf{E}\big[\big(\tau_0^d- (N-K)\big)^+\big]}{\mathbf{E}[\tau_0^d]} ,
 \end{equation*}
 which vanishes as $N\to\infty$.
 \qed

 Before proving that the second term in (\ref{eqn4.08}) also can be ignored, we need to provide the following important ancillary result, which provides a uniformly estimates for the probability $\mathbf{P}_x\{ \tau_0^d\geq n  \}$  starting from different state $x$.

 \begin{proposition}\label{Thm407}
 For any $\epsilon>0 $ satisfying that $ \epsilon\cdot \mathbf{E}[\tau_{-1}^{d}]<1$, there exists a constant $C_\epsilon>0$ such that for any $n\geq 1$ and $x\in[0,\epsilon n]$,
 \begin{equation*}
 \mathbf{P}_x\{ \tau_0^d\geq n  \} \leq C_\epsilon(1+x) \cdot \mathbf{P}\{ \xi_1> an \}.
 \end{equation*}
 \end{proposition}
 \proof  Let  $\{\tau_{-1}^{d,i}\}_{i\geq 1}$ be a sequence of i.i.d copies of $\tau_{-1}^d$. From the strong Markov property of $\xi$, for any $\theta>1$ with $\theta\cdot \epsilon\cdot \mathbf{E}[\tau_{-1}^{d}]<1$ we have as $n\to\infty$,
 \begin{eqnarray*}
 \mathbf{P}_x\{ \tau_0^d\geq n  \} \leq \mathbf{P}\Big\{\sum_{i=1}^{[x]+1} \tau_{-1}^{d,i}\geq n  \Big\}\ar= \ar\mathbf{P}\Big\{\sum_{i=1}^{[x]+1} \big(\tau_{-1}^{d,i}-\theta  \mathbf{E}[\tau_{-1}^{d}]\big)  \geq n \big(1- \theta  \mathbf{E}[\tau_{-1}^{d}]\cdot (1+[x])/n\big) \Big\}\cr
 \ar\leq\ar \mathbf{P}\Big\{\sum_{i=1}^{[x]+1} \big(\tau_{-1}^{d,i}-\theta  \mathbf{E}[\tau_{-1}^{d}]\big)  \geq n \big(1- \theta  \mathbf{E}[\tau_{-1}^{d}]\cdot [\epsilon n]/n\big) \Big\}\cr
 \ar\sim \ar \mathbf{P}\Big\{\sum_{i=1}^{[x]+1} \big(\tau_{-1}^{d,i}-\theta  \mathbf{E}[\tau_{-1}^{d}]\big)  \geq n \big(1- \epsilon\cdot \theta  \mathbf{E}[\tau_{-1}^{d}]\big) \Big\}.
 \end{eqnarray*}
 From Theorem 2 in \cite{DenisovFossKorshunov2010},  there exists a constant $C>0$ such that for any $x\geq 0$,
 \begin{eqnarray*}
 \lefteqn{\mathbf{P}\Big\{\sum_{i=1}^{[x]+1} \big(\tau_{-1}^{d,i}-\theta  \mathbf{E}[\tau_{-1}^{d}]\big)  \geq n  \big(1- \epsilon\cdot \theta  \mathbf{E}[\tau_{-1}^{d}]\big) \Big\}	}\ar\ar\cr
 \ar\leq\ar C(x+1) \mathbf{P}\big\{ \tau_{-1}^{d}-\theta  \mathbf{E}[\tau_{-1}^{d}]\geq n \big(1- \epsilon\cdot \theta  \mathbf{E}[\tau_{-1}^{d}]\big)  \big\}
 \leq C(x+1) \mathbf{P}\big\{ \tau_{-1}^{d} \geq n \big(1- \epsilon\cdot \theta  \mathbf{E}[\tau_{-1}^{d}]\big)  \big\}.
\end{eqnarray*}
From  (\ref{eqn3.21}), we have
 \begin{eqnarray*}
\mathbf{P}\big\{ \tau_{-1}^{d} \geq n \big(1- \epsilon\cdot \theta  \mathbf{E}[\tau_{-1}^{d}]\big)  \big\} \ar\sim \ar C \mathbf{P}\big\{ \xi_1> n\cdot  \big(1- \epsilon\cdot \theta  \mathbf{E}[\tau_{-1}^{d}]\big)  \big\}
\end{eqnarray*}
and hence
 \begin{eqnarray*}
 \mathbf{P}\Big\{\sum_{i=1}^{[x]+1} \big(\tau_{-1}^{d,i}-\theta  \mathbf{E}[\tau_{-1}^{d}]\big)  \geq n  \big(1- \epsilon\cdot \theta  \mathbf{E}[\tau_{-1}^{d}]\big) \Big\}
 \ar\leq \ar  C \big(1- \epsilon\cdot \theta  \mathbf{E}[\tau_{-1}^{d}]\big)^{-\alpha}\cdot (x+1)\mathbf{P}\{ \xi_1> an \}.
 \end{eqnarray*}
 Here we have gotten the desired result.
 \qed

 The following lemma considers the second term on the right side of (\ref{eqn4.08}) with the observation that sample paths with late arrival local infimum usually result in the fast increasing of $A_t(\xi)$ and hence their contribution to the expectation (\ref{eq1.2}) decreases fast.

 \begin{lemma}\label{Thm406}
 There exists a constant $C>0$ depending on the quantity $K_0$ in Condition~\ref{ConF1} such that for any $t\geq 0$,
 \begin{equation}\label{eqn4.10}
 \mathbf{E}[F(A_t(\xi))]\leq C\mathbf{P}\{\xi_1>at\}.
 \end{equation}
 Moreover, we also have
 \begin{equation}\label{eqn4.11}
 \lim_{K\rightarrow \infty}\limsup_{n\rightarrow\infty}\frac{\mathbf{E}[F(A_{n+1}(\xi));\, g^d_n\geq K]}{\mathbf{P}\{\xi_1>an\}}=0.
 \end{equation}
 \end{lemma}
 \proof
 We first prove the second result. By the monotonicity of $F$, we have
 \begin{eqnarray}\label{eqn4.12}
  \mathbf{E}[F(A_{n+1}(\xi)),\, g^d_n\geq K] \ar\leq \ar \sum_{k=K}^n\mathbf{E}\Big[F\Big(e^{-\xi_k}\int_0^1e^{-(\xi_{k+s}-\xi_k)}ds\Big); \, g^d_n=k\Big] \cr
 \ar\leq\ar \sum_{k=K}^{n}\mathbf{E}\Big[\mathbf{E}\Big[F\big(e^{-\xi_k}A_1(\xi')\big); \inf_{i=1,\cdots,n-k}\xi'_{i}> 0\,\Big|\,  \sigma(\xi)\Big] ; \, g^d_k=k\Big],\qquad
 \end{eqnarray}
 where $\sigma(\xi):= \sigma(\xi_t:t\geq 0)$ and  $\hat\xi$ is an independent copy of $\xi$.
 For any $\epsilon>0 $ satisfying that $ \epsilon\cdot \mathbf{E}[\tau_{-1}^{d}]<1$, we have
 \begin{eqnarray}\label{eqn4.13}
 \lefteqn{\mathbf{E}\Big[F\big(e^{-\xi_k}A_1(\xi')\big); \,\inf_{k=1,\cdots, n-k}\xi'_{i}> 0\,\Big|\,  \sigma(\xi)\Big]}\ar\ar\cr
 \ar=\ar
 \mathbf{E}\Big[F\big(e^{-\xi_k}A_1(\xi')\big) ; \,
 0< \xi'_1\leq \epsilon (n-k) ,\, \inf_{k=2,\cdots, n-k}\xi'_{i}> 0 \,\Big|\,  \sigma(\xi)\Big]\cr
 \ar\ar +\mathbf{E}\Big[F\big(e^{-\xi_k}A_1(\xi')\big) ; \,
 \xi'_1>\epsilon (n-k) ,\, \inf_{k=2,\cdots, n-k}\xi'_{i}> 0 \,\Big|\,  \sigma(\xi)\Big].
 \end{eqnarray}
 Let $\mathcal{N}_1^{'x}:= \#\{ t\in[0,1]: \Delta \xi'_t>x \}$ for any $x>0$.
  From the first result in (\ref{eqn4.07}), we have $\xi'_1>\epsilon (n-k) $ for large $n$ if and only if $\mathcal{N}_1^{'\epsilon (n-k)}=1 $.  From this and the boundedness of $F$, we see that the second term on the right side of above equation can be bounded by
  \begin{eqnarray*}
  	\mathbf{E}\Big[F\big(e^{-\xi_k}A_1(\xi')\big) ; \,
  	\xi'_1>\epsilon (n-k) \,\Big|\, \sigma(\xi)\Big]
  	\ar\sim\ar  \mathbf{E}\Big[F\big(e^{-\xi_k}A_1(\xi')\big)  ; \,
  	\mathcal{N}_1^{'\epsilon (n-k)}=1\,\Big|\, \sigma(\xi) \Big] .
  \end{eqnarray*}
  From Condition~\ref{ConF1}, we also have
  \begin{eqnarray*}
 \mathbf{E}\Big[F\big(e^{-\xi_k}A_1(\xi')\big)  ; \,
  	\mathcal{N}_1^{'\epsilon (n-k)}=1 \,\Big|\, \sigma(\xi)\Big]
  	\ar\leq\ar  K_0\cdot e^{\beta\xi_k}\cdot \mathbf{E}\Big[\big(A_1(\xi')\big)^{-\beta} ; \,	\mathcal{N}_1^{'\epsilon (n-k)}=1 \Big]\cr
  	\ar=\ar K_0\cdot e^{\beta\xi_k} \cdot \mathbf{E}\Big[\big(A_1(\xi)\big)^{-\beta} ; \,	\mathcal{N}_1^{\epsilon (n-k)}=1 \Big].
  \end{eqnarray*}
  Here conditoned on $\mathcal{N}_1^{\epsilon (n-k)}=1$, we have $\mathcal{J}^{\epsilon (n-k)}\in[0,1]$ a.s.  and
    \begin{eqnarray*}
   \mathbf{E}\Big[\big(A_1(\xi)\big)^{-\beta} ; \,	\mathcal{N}_1^{\epsilon (n-k)}=1 \Big]
   \ar\leq\ar
    \mathbf{E}\Big[\big(A_{\mathcal{J}^{\epsilon (n-k)}}(\xi^{\epsilon (n-k)})\big)^{-\beta} \,\Big|\, \mathcal{N}_1^{\epsilon (n-k)}=1 \Big] \cdot\mathbf{P}\{ \mathcal{N}_1^{\epsilon (n-k)}=1  \}.
  \end{eqnarray*}
 From the independence between $\mathcal{J}^{\epsilon (n-k)}$ and $\xi^{\epsilon (n-k)}$, there exists a uniformly distributed random variable $\mathbf{U}$ on $[0,1]$ independent of $\xi$ such that
 \begin{eqnarray*}
  \mathbf{E}\Big[\big(A_{\mathcal{J}^{\epsilon (n-k)}}(\xi^{\epsilon (n-k)})\big)^{-\beta} \,\Big|\, \mathcal{N}_1^{\epsilon (n-k)}=1 \Big]
  \ar=\ar  \mathbf{E}\Big[\big(A_{\mathbf{U}}(\xi^{\epsilon (n-k)})\big)^{-\beta}   \Big]  \leq \int_0^1\mathbf{E}\Big[\big(A_s(\xi )\big)^{-\beta}   \Big] ds,
 \end{eqnarray*}
 which is finite because of $\mathbf{E}\big[\big(A_s(\xi )\big)^{-\beta}   \big]\sim s^{-\beta}$ as $s\to 0+$; see Theorem 2.18(1) in \cite{PatieSavov2018}.
 Putting all result above together, we have
 \begin{eqnarray}\label{eqn4.13.01}
  \mathbf{E}\Big[F\big(e^{-\xi_k}A_1(\xi')\big) ;\, \xi'_1>\epsilon (n-k) \,\Big|\, \sigma(\xi) \Big] \ar\leq\ar C e^{\beta\xi_k} \cdot\mathbf{P}\{ \mathcal{N}_1^{\epsilon (n-k)}=1  \}\cr
  \ar\leq\ar C(a/\epsilon)^\alpha \cdot e^{\beta\xi_k} \cdot \mathbf{P}\{ \xi_1> a (n-k)  \}.
 \end{eqnarray}
 From the Markov property of $\xi'$, we see that the first term on the right side of (\ref{eqn4.13}) equals to
 \begin{eqnarray*}
  \lefteqn{\mathbf{E}\Big[F\big(e^{-\xi_k}A_1(\xi')\big)\cdot \mathbf{P}\Big\{\inf_{k=2,\cdots, n-k}\xi'_{i}\geq 0\,\Big|\, \xi'_1\Big\} ; \, 0<\xi'_1\leq \epsilon (n-k)  \,\Big|\, \sigma(\xi)\Big] }\ar\ar\cr
  \ar\leq\ar  C \mathbf{E}\Big[F\big(e^{-\xi_k}A_1(\xi')\big)\cdot
  (\xi'_1\vee 0+1) \,\Big|\, \sigma(\xi) \Big]\cdot  \mathbf{P} \big\{\xi_1> a(n-k-1)\big\}.
 \end{eqnarray*}
Here the inequality above follows from Proposition~\ref{Thm407}.
From Lemma~\ref{Thm403}, there exists a constant $\lambda_0\in(0,\beta)$ such that
\begin{eqnarray*}
 \mathbf{E}\Big[F\big(e^{-\xi_k}A_1(\xi')\big)\cdot
	(\xi'_1\vee 0+1) \,\Big|\,  \sigma(\xi)\Big]
	\ar\leq\ar Ce^{\lambda_0 \xi_k}
\end{eqnarray*}
and  the first term on the right side of (\ref{eqn4.13}) can be bounded by
 \begin{eqnarray*}
	 Ce^{\lambda_0 \xi_k} \cdot \mathbf{P} \big\{\xi_1> a(n-k-1) \big\} .
\end{eqnarray*}
 Taking this and (\ref{eqn4.13.01}) back into (\ref{eqn4.13}), we have
  \begin{eqnarray*}
 \lefteqn{\mathbf{E}\Big[\mathbf{E}\Big[F\big(e^{-\xi_k}A_1(\xi')\big) ; \, \inf_{i=1,\cdots,n-k}\xi'_{i}\geq 0\,\Big|\,  \sigma(\xi)\Big] ; \, g^d_k=k\Big]}\ar\ar\cr
 \ar\ar\cr
 \ar\leq\ar
 C \mathbf{E}\big[e^{\lambda_0  \xi_k}; g^d_k=k\big]\cdot \mathbf{P} \{\xi_1>   a(n-k-1)\}.
 \end{eqnarray*}
 By the duality of $\xi$ and  (\ref{eqn3.44}),   we have $\{\xi_{k-i}- \xi_k:i=0,\cdots,k\}\overset{\rm d}= \{ -\xi_i:i=0,\cdots,k  \}$ and hence   as $k\to\infty$,
 \begin{eqnarray*}
 	\mathbf{E}\big[e^{\lambda_0  \xi_k}; g^d_k=k\big]	\ar=\ar 	\mathbf{E}\Big[e^{\lambda_0  \xi_k}; \inf_{i=0,\cdots, k-1}\xi_i \geq \xi_k\Big]\cr
  	\ar=\ar 	\mathbf{E}\Big[e^{-\lambda_0  (\xi_0-\xi_k)};\inf_{i=0,\cdots k}(\xi_{k-i}- \xi_k)\geq  0\Big]\cr
 	\ar=\ar	\mathbf{E}\Big[e^{\lambda_0 \xi_k};\tau_0^{d+}> k\Big]\sim  \frac{C}{ak} \cdot \mathbf{P} \{\xi_1> ak\}.
 	\end{eqnarray*}
 From all estimate above and Proposition~\ref{Thm111}, we have for large $n$,
 \begin{eqnarray*}
  \mathbf{E}[F(A_{n+1}(\xi)) ; \, g^d_n\geq K]
 \ar\leq\ar C\cdot  \sum_{k=K}^{n} \frac{1}{k} \mathbf{P}\{\xi_1\geq ak\} \cdot \mathbf{P} \{\xi_1> a(n-k-1)\} \cr
 \ar\leq \ar \frac{C}{K}\cdot  \sum_{k=0}^{n}  \mathbf{P}\{\xi_1\geq ak\} \cdot \mathbf{P} \{\xi_1> a(n-k-1)\}  \cr
 \ar\sim \ar \frac{C}{K} \sum_{k=0}^{\infty} \mathbf{P}\{\xi_1\geq ak\} \cdot \mathbf{P} \{\xi_1> an\}.
 \end{eqnarray*}
 Taking this back into (\ref{eqn4.12}),  we have
 \begin{eqnarray}\label{eqn4.12.1}
 \limsup_{n\to\infty}\frac{ \mathbf{E}[F(A_{n+1}(\xi));g^d_n\geq K]}{\mathbf{P} \{\xi_1> an\}}\leq \frac{C}{K},
 \end{eqnarray}
 which vanishes as $K\to\infty$. Here we have got the second result. We now start to prove the first one. Actually, we have
 \begin{eqnarray*}
  \limsup_{n\to\infty} \frac{ \mathbf{E}[F(A_{n+1}(\xi))]}{\mathbf{P} \{\xi_1> an\}}
  \ar=\ar  \limsup_{n\to\infty}\frac{ \mathbf{E}[F(A_{n+1}(\xi));g^d_n\geq 1]}{\mathbf{P} \{\xi_1> an\}}+ \limsup_{n\to\infty}\frac{ \mathbf{E}[F(A_{n+1}(\xi));g^d_n=0]}{\mathbf{P} \{\xi_1> an\}}\cr
  \ar\leq\ar C + C  \limsup_{n\to\infty} \frac{ \mathbf{P}\{ \tau_0^d\geq n  \}  }{ \mathbf{P} \{\xi_1> an\}  } <\infty.
 \end{eqnarray*}
Here the first inequality follows from (\ref{eqn4.12.1}) and the boundedness of $F$, and the second inequality follows directly from (\ref{eqn3.21}).
This together with the non-increasing of $F$ immediately induces the  first desired result.
 \qed

 \subsection{Contribution of sample paths with early arrival of large jump}

 We now start to analyze the contribution of sample paths with early arrival large jump to the expectation (\ref{eq1.2}).
 In this case, we observe that the effect of their partial paths before the large jump on $A_{t}(\xi)$ is slight and hence it is the key step to clarify the increasing rate of $A_{t}(\xi)$ after the large jump.
 As we  have showed in Theorem~\ref{Thm108}, the sample paths stay in the set far away from origin for a log time after the large jump. This suggests us that $A_t(\xi)$ can be well approximated by $A_T(\xi)$ with $T$ larger than the arrival time of the first large jump; see the following lemma.

 \begin{lemma}\label{Thm409}
 If $F$ is globally Lipschitz continuous on $(0,\infty)$, we have
 \begin{equation}\label{eqn4.14}
 \lim_{T\rightarrow \infty}\limsup_{t\rightarrow \infty} \mathbf{E}\big[|F(e^{- \mathbf{J}}A_{T}(\xi))-F(e^{-\mathbf{J}}A_{t}(\xi))|\, \big|\, \mathbf{J}> at\big]= 0.
 \end{equation}
 \end{lemma}
 \proof
 Since $F$ is bounded and globally Lipschitz continuous, there exist two constants $K_1,K_2>0$ such that for any $b>a$
 \begin{equation*}
 \mathbf{E}\big[|F(e^{-\mathbf{J}}A_{T}(\xi))-F(e^{-\mathbf{J}}A_{t}(\xi))| \, \big|\, \mathbf{J}> at\big]
 \leq \mathbf{E}\Big[\Big( K_1 e^{-\mathbf{J}}\int_T^te^{-\xi_s}ds\Big)\wedge K_2 \,\Big|\,\mathbf{J}> at\Big],
 \end{equation*}
 which can be bounded by $\varepsilon_1(b,t)+\varepsilon_2(b,T,t)$ with
 \begin{equation*}
 \varepsilon_1(b,t):= \frac{K_2\mathbf{P}\{\mathbf{J}\in (at,bt]\}}{\mathbf{P}\{\mathbf{J}>at\}}
 \quad\mbox{and}\quad
 \varepsilon_2 (b,T,t)
 :=\frac{\mathbf{E}[(K_1 e^{-\mathbf{J}}\int_T^te^{-\xi_s}ds)\wedge K_2 ; \mathbf{J}> bt]}{\mathbf{P}\{\mathbf{J}>at\}}.
 \end{equation*}
 From the definition of $\mathbf{J}$, we have $ \varepsilon_1(b,t) \to K_2[1-(a/b)^\alpha]$ as $t\to\infty$.
 Moreover,
 \begin{equation*}
 \varepsilon_2 (b,T,t)
 \leq \mathbf{E}\Big[\Big(K_1 \int_T^te^{- (\xi_s+bs)}ds\Big)\wedge K_2\Big]\cdot \frac{\mathbf{P}\{\mathbf{J}> bt\}}{\mathbf{P}\{\mathbf{J}> at\}}
 \leq \mathbf{E}\Big[\Big(K_1 \int_T^\infty e^{- (\xi_s+bs)}ds\Big)\wedge K_2\Big],
 \end{equation*}
 which vanishes as $T\to\infty$ because of $\mathbf{E}[\xi_1]+b>0$  and $\int_0^\infty e^{- (\xi_s+bs)}ds<\infty$ a.s.
 Putting these estimates together, we see that  (\ref{eqn4.14}) follows as $b\to a+$.
 \qed

 \begin{lemma}\label{Thm408}
 For any $T>0$, there exists a constant $C_{F,T}>0$ such that
 \begin{equation}\label{eqn4.15}
 \lim_{t\to \infty}\mathbf{E}\big[F(e^{- \mathbf{J}}A_{T}(\xi))\,\big|\,\mathbf{J}> at\big]= C_{F,T}.
 \end{equation}
 Moreover, $C_{F,T}$ decreases to $C_F>0$  as  $T\to \infty$  and
 \begin{equation}\label{eqn4.16}
 \lim_{t\rightarrow \infty}\mathbf{E}\big[F(e^{- \mathbf{J}}A_t(\xi))\,\big|\,\mathbf{J}> at\big]= C_F.
 \end{equation}
 \end{lemma}
 \proof
 Since $F$ is  non-increasing, we have $\mathbf{E}[F(e^{- \mathbf{J}}A_{T}(\xi))\,|\,\mathbf{J}\geq at]$ is non-decreasing in $t$. Indeed, we have
 \begin{eqnarray*}
 	\mathbf{E}[F(e^{- \mathbf{J}}A_{T}(\xi))\,|\,\mathbf{J}\geq at] = \int_0^\infty 	\mathbf{E}[F(e^{- \mathbf{J}}\cdot y\,|\,\mathbf{J}\geq at] \mathbf{P}\{  A_{T}(\xi)\in dy \}
 	\end{eqnarray*}
 and for any $s<t$ and $y\geq 0$,
  \begin{eqnarray*}
 \lefteqn{ \mathbf{E}[F(e^{- \mathbf{J}}\cdot y\,|\,\mathbf{J}\geq as]
  =	\int_{as}^{at} F(e^{-x} y) \frac{\mathbf{P}\{ \mathbf{J} \in dx \}}{\mathbf{P}\{ \mathbf{J}\geq as\}} + 	\frac{\mathbf{P}\{ \mathbf{J}\geq at\}}{\mathbf{P}\{ \mathbf{J}\geq as\}}\int_{at}^\infty F(e^{-x} y) \frac{\mathbf{P}\{ \mathbf{J} \in dx \}}{\mathbf{P}\{ \mathbf{J}\geq at\} }   }\qquad \qquad\ar\ar\cr
  	\ar=\ar \int_{at}^\infty F(e^{-x} y) \frac{\mathbf{P}\{ \mathbf{J} \in dx \}}{\mathbf{P}\{ \mathbf{J}\geq at\} }\cr
  \ar\ar	+ \frac{\mathbf{P}\{ \mathbf{J}\in [as,at)\}}{\mathbf{P}\{ \mathbf{J}\geq as\}}\Big(\int_{as}^{at}  \frac{F(e^{-x} y)\mathbf{P}\{ \mathbf{J} \in dx \}}{\mathbf{P}\{ \mathbf{J}\in [as,at)\}}  -  \int_{at}^\infty \frac{F(e^{-x} y) \mathbf{P}\{ \mathbf{J} \in dx \}}{\mathbf{P}\{ \mathbf{J}\geq at\} }\Big)\cr
  	\ar=\ar   \mathbf{E}[F(e^{- \mathbf{J}}\cdot y\,|\,\mathbf{J}\geq at]  + \frac{\mathbf{P}\{ \mathbf{J}\in [as,at)\}}{\mathbf{P}\{ \mathbf{J}\geq as\}}
  	\big(  \mathbf{E}[F(e^{- \eta_1}\cdot y)]- \mathbf{E}[F(e^{- \eta_2}\cdot y)]
  	\big),
 \end{eqnarray*}
where $\eta_1$ and $\eta_2$ are two random variables satisfying that
 \begin{eqnarray*}
 \mathbf{P}\{ \eta_1\in dx  \} =	\mathbf{1}_{x\in [as,at)} \frac{\mathbf{P}\{ \mathbf{J} \in dx \}}{\mathbf{P}\{ \mathbf{J}\in [as,at)\}}
 \quad \mbox{and}\quad
  \mathbf{P}\{ \eta_2\in dx  \} =	\mathbf{1}_{x\geq at} \frac{\mathbf{P}\{ \mathbf{J} \in dx \}}{\mathbf{P}\{ \mathbf{J}\geq at\}} .
\end{eqnarray*}
Since $F$  non-increasing, we see that $ \mathbf{E}[F(e^{- \eta_1}\cdot y)]\leq  \mathbf{E}[F(e^{- \eta_2}\cdot y)]  $ and hence $\mathbf{E}[F(e^{- \mathbf{J}}A_{T}(\xi))\,|\,\mathbf{J}\geq at]$ is non-decreasing in $t$.
Moreover, since $F$ is bounded, we also have $\mathbf{E}[F(e^{- \mathbf{J}}A_{T}(\xi))\,|\,\mathbf{J}\geq at]$ is uniformly bounded and hence the limit (\ref{eqn4.15}) holds.
 We now start to prove the second result.  Notice that $C_{F,T}$ is non-increasing in $T$ and hence $C_{F,T} \to C_{F}\in[0,\infty)$ as $T\rightarrow \infty$.
  From Lemma~\ref{Thm409}, we have
  \begin{eqnarray*}
 	 \lim_{t\rightarrow \infty}\mathbf{E}\big[F(e^{-\mathbf{J}}A_{t}(\xi))\,\big|\,\mathbf{J}> at\big]
 	= \lim_{T\to\infty}  \lim_{t\rightarrow \infty}	\mathbf{E}\big[F(e^{-\mathbf{J}}A_T(\xi))\,\big|\,\mathbf{J}> at\big]=\lim_{T\to\infty} C_{F,T}= C_F.
 \end{eqnarray*}
 We now show that $C_F>0$.
 For any $b>a$, let $\tilde{\xi}_t=\xi_t+ bt$, which drifts to $\infty $. We have
 \begin{eqnarray*}
 \mathbf{E}\big[F(e^{-\mathbf{J}}A_{t}(\xi))\,\big|\,\mathbf{J}> at\big]\ar\geq \ar \frac{\mathbf{E}[F(e^{-\mathbf{J}}A_{t}(\xi)),\mathbf{J}> bt]}{\mathbf{P}\{\mathbf{J}> at\}}\cr
 \ar\geq\ar  \mathbf{E}\big[F(e^{- bt}A_{t}(\xi))\big] \frac{\mathbf{P}\{\mathbf{J}> bt\}}{\mathbf{P}\{\mathbf{J}> at\}}
  \geq \mathbf{E}\big[F(A_{t}(\tilde\xi))\big] \frac{\mathbf{P}\{\mathbf{J}> bt\}}{\mathbf{P}\{\mathbf{J}> at\}}.
 \end{eqnarray*}
 Since $ \mathbf{E}[F(A_\infty(\tilde\xi))]\in(0,\infty) $,  there exists   $\underline{C}>0$ such that for any $T>0$,
 \begin{equation}
 \liminf_{t\to\infty}  \mathbf{E}\big[F(e^{- \mathbf{J}}A_{T}(\xi))\,\big|\,\mathbf{J}> at\big]\geq \underline{C},
 \end{equation}
 which immediately induces that $C_F>0$.
 \qed

 \subsection{Proof for Theorem~\ref{T1}}

 We first consider the special case in which $F(x)$ is globally Lipschitz continuous. We first have
 \begin{equation*}
 \lim_{t\to\infty}\frac{\mathbf{E}[F(A_t(\xi))]}{\mathbf{P}\{\xi_1>at\}}
 = \lim_{N\to\infty}\lim_{t\to\infty}\frac{\mathbf{E}[F(A_t(\xi));\mathcal{J}^{at}\leq N]}{\mathbf{P}\{\xi_1>at\}} +\lim_{N\to\infty}\lim_{t\to\infty}\frac{\mathbf{E}[F(A_t(\xi));\mathcal{J}^{at}>N ]}{\mathbf{P}\{\xi_1>at\}}.
 \end{equation*}
 From Lemma~\ref{Thm411} and \ref{Thm406}, the second limit on the right side of above equation  equals to $0$.
 For large $t>N$,
 \begin{equation*}
  \mathbf{E}\big[F(A_t(\xi));\mathcal{J}^{at}\leq N\big]
  =\mathbf{E}\Big[F\Big(\int_{(0,\mathcal{J}^{at})}e^{-\xi^{at}_r}dr+e^{-\xi^{at}_{\mathcal{J}^{at}-}-\Delta\xi_{\mathcal{J}^{at}}}\int_{\mathcal{J}^{at}}^te^{-(\xi_r-\xi_{\mathcal{J}^{at}})}dr\Big);\mathcal{J}^{at}\leq N\Big].
 \end{equation*}
 From (\ref{eqn3.04}) and the independence between $\mathcal{J}^{at}$ and  $\Delta\xi_{\mathcal{J}^{at}}$, we have
 \begin{equation*}
 \mathbf{P}\{  \mathcal{J}^{at} \in ds, \Delta\xi_{\mathcal{J}^{at}}\in dy  \} =\mathbf{1}_{\{y\geq at\}}  e^{-\bar\nu(at)\cdot s}ds \nu(dy) .
 \end{equation*}
 From the independent increments of $\xi$, we have
 \begin{eqnarray*}
  \lefteqn{ \mathbf{E}\big[F(A_t(\xi));\mathcal{J}^{at}\leq N\big]} \ar\ar\cr
  \ar\ar\cr
   \ar=\ar \int_0^N \int_{at}^\infty \mathbf{E}\Big[F\Big(\int_{(0,s)}e^{-\xi^{at}_r}dr+e^{-\xi^{at}_{s-}-y}\int_0^{t-s}e^{-\hat\xi_r}dr\Big) \Big]\mathbf{P}\{  \mathcal{J}^{at} \in ds, \Delta\xi_\mathcal{J}^{at}\in dy  \}\cr
   \ar=\ar \int_0^N \bar\nu(at)e^{-\bar\nu(at)\cdot s}ds\int_{at}^\infty \mathbf{E}\Big[  F\big(A_{s-}(\xi^{at})+e^{-\xi^{at}_{s-}-y}A_{t-s}(\hat\xi)\big) \Big]  \frac{\nu(dy) }{\bar\nu(at)}\cr
   \ar=\ar \bar\nu(at) \int_0^N e^{-\bar\nu(at)\cdot s} \cdot  \mathbf{E}\Big[ \mathbf{E}\big[ F\big(A_{s-}(\xi^{at})+e^{-\xi^{at}_{s-}-\mathbf{J}}A_{t-s}(\hat\xi)\big) \big|\sigma(\xi),  \mathbf{J}\geq at\big] \Big]ds .
  \end{eqnarray*}
   By the dominated convergence theorem and (\ref{eqn3.01})-(\ref{eqn3.02}), we have
 \begin{eqnarray*}
 \lim_{t\to\infty}\frac{  \mathbf{E}\big[F(A_t(\xi));\mathcal{J}^{at}\leq N\big]}{\mathbf{P}\{\xi_1>at\}}\ar=\ar   \int_0^N   \mathbf{E}\Big[   \lim_{t\to\infty}\mathbf{E}\big[ F\big(A_{s-}(\xi^{at})+e^{-\xi^{at}_{s-}-\mathbf{J}}A_{t-s}(\hat\xi)\big) \big| \sigma(\xi), \mathbf{J}>at\big] \Big] ds .
 \end{eqnarray*}
 From Lemma~\ref{Thm408} and the fact that $\sup_{r\in[0,s)}|\xi^{at}_r-\xi_r|\to 0$ a.s.,  the limit above exists and hence by  the stochastic continuity of $A_\cdot(\xi)$ and $\xi_\cdot$,
 \begin{eqnarray*}
  \lim_{t\to\infty}\frac{\mathbf{E}[F(A_t(\xi))]}{\mathbf{P}\{\xi_1>at\}}\ar=\ar \int_0^\infty \mathbf{E}\big[ C_F(\xi,s-)\big]ds = \int_0^\infty \mathbf{E}\big[ C_F(\xi,s)\big]ds,
 \end{eqnarray*}
 which  is finite; see (\ref{eqn4.10}).
 We now prove this theorem for general $F$.
 For $n\ge 1$, let $F_n(y) = F(1/n)1_{\{y\le 1/n\}} + F(y)1_{\{y>1/n\}}$, which is globally Lipschitz and $F_n(y)\to F(y)$ increasingly.
 From the result above, we have
 \begin{eqnarray*}
 \lim_{t\to\infty}\frac{\mathbf{E}[F_n(A_t(\xi))]}{\mathbf{P}\{\xi_1>at\}}
 = \int_0^\infty\mathbf{E}[C_{F_n}(\xi,s)]ds
 \end{eqnarray*}
 with $\mathbf{E}[C_{F_n}(\cdot)]$ is non-decreasing in $n$.
 From the monotone convergence theorem and (\ref{eqn4.10}), we have as $n\to\infty$,
 \begin{eqnarray*}
  \int_0^\infty\mathbf{E}[C_{F_n}(\xi,s)]ds\to \int_0^\infty \mathbf{E}[C_F(\xi, s)]ds<\infty.
 \end{eqnarray*}
 Let $G_n(y) = F(y)-F_n(y)$.  From Condition~\ref{ConF1} and \ref{ConF2}, it is easy to see  for any $n\geq 1$ and $x\geq 0$ have $G_n(x)/G_n(0)\leq 1 \wedge  x^{-\beta}.$
 From Lemma~\ref{Thm406}, we can prove that there exists a constant $C>0$ such that for any $n\geq 1$ and $t\geq 0$,
 \begin{eqnarray*}
 \mathbf{E}[G_n(A_t(\xi))/G_n(0)]\leq C\mathbf{P}\{\xi_1>at\}
 \quad \mbox{and}\quad
 \limsup_{t\to\infty}\frac{\mathbf{E}[G_n(A_t(\xi))]}{\mathbf{P}\{\xi_1>at\}}\leq CG_n(0),
 \end{eqnarray*}
 which goes to $0$ as $n\to\infty$. Putting all results together, we have
 \begin{eqnarray*}
 \lim_{t\to\infty}\frac{\mathbf{E}[F(A_t(\xi))]}{\mathbf{P}\{\xi_1>at\}} =\lim_{n\to\infty}\lim_{t\to\infty}\frac{\mathbf{E}[G_n(A_t(\xi))]}{\mathbf{P}\{\xi_1>at\}}
 +\lim_{n\to\infty}\int_0^\infty \mathbf{E}[C_{F_n}(\xi,s)]ds= \int_0^\infty \mathbf{E}[C_F(\xi,s)]ds.
 \end{eqnarray*}
 Here we have finished the proof.


 \section{Asymptotic results for CBRE-processes}
 \setcounter{equation}{0}

 In this section, we apply the results in the last section to study the asymptotic behavior of survival probabilities of continuous-state branching processes in L\'evy random environment.
 Let $\{Z^\gamma_t: t\ge 0\}$ be a spectrally positive $(\gamma+1)$-stable process with $0<\gamma\leq 1$ and $\{Z^e_t: t\ge 0\}$ be a L\'{e}vy process with no jump less than $-1$.
 When $\gamma=1$, we think of $\{Z^\gamma(t): t\ge 0\}$ as a Brownian motion. When $0< \gamma< 1$, we assume $\{Z^\gamma(t): t\ge 0\}$ has L\'{e}vy measure $m(dz) =\frac{\gamma1_{\{z>0\}}dz}{\Gamma(1-\gamma)z^{2+\gamma}}.$
 Associated to the L\'evy processes $\xi$ defined by (\ref{eqn2.05}), we may assume the L\'evy process $Z^e$ admits the following L\'evy-It\^o decomposition:
 \begin{eqnarray}\label{eqn5.02}
 Z^e_t =  a_0 t + \sigma B_t + \int_0^t \int_{[-1,1]}(e^u-1)\tilde{N}(ds,du) + \int_0^t \int_{[-1,1]^c} (e^u-1) N(ds,du),
 \end{eqnarray}
 where $[-1,1]^c = \mathbb{R}\setminus [-1,1]$, $\tilde{N}(ds,du) = N(ds,du) - ds\nu(du)$ and
 \begin{eqnarray*}
 a_0 = -a+\frac{\sigma^2}{2} + \int_{[-1,1]}(e^z-1-z) \nu(dz) - \int_{[-1,1]^c}z\nu(dz).
 \end{eqnarray*}
 Then $Z^e$ has no jump smaller than $-1$.
 Clearly, the two processes $Z^e$ and $\xi$ generate the same filtration.
 Let $c\ge 0$ be another constant. Given the initial value $x\ge 0$, by Theorem~6.2 in \cite{FuLi2010}, there exists a unique positive strong solution $\{X_t(x):t\geq 0\}$ to (\ref{eq1.3}). The solution is called a \textit{continuous-state branching process in random environment (CBRE-process)} with \textit{stable branching mechanism}. Here the random environment is modeled by the L\'evy process $Z^e$.
 The reader may refer to \cite{HeLiXu2018} and \cite{PalauPardo2018} for discussions of more general CBRE-processes.

 Let $\mathbf{P}^\xi$ denote the conditional law given $Z^e$ or $\xi$.
 Let $Z_t(x) = X_t(x)\exp\{-\xi_t\}$ for any $t\geq 0$.
 For any $  \lambda\ge 0 $ and $t\ge r\ge 0,$ by Theorem~1 in \cite{BansayePardoSmadi2013} or Theorem~3.4 in \cite{HeLiXu2018}  we have
 \begin{eqnarray}\label{eqn5.03}
 \mathbf{E}^\xi [e^{-\lambda Z_t(x)}|\mathscr{F}_r] = \exp\{-Z_r(x)u^\xi_{r,t}(\lambda)\},
 \end{eqnarray}
 where $r\mapsto u^\xi_{r,t}(\lambda)$ is the solution to
 \begin{eqnarray*}
 \frac{d}{dr}u^\xi_{r,t}(\lambda) = ce^{-\gamma \xi(r)} u^\xi_{r,t}(\lambda)^{\gamma+1},
 \qquad
 u^\xi_{t,t}(\lambda) = \lambda.
 \end{eqnarray*}
 By solving the above equation, we get
 \begin{eqnarray}\label{eqn5.04}
 u^\xi_{r,t}(\lambda)
 =\Big(c\gamma\int_r^t e^{-\gamma \xi(s)}ds + \lambda^{-\gamma}\Big)^{-1/\gamma};
 \end{eqnarray}
 see the proof of Proposition~4 in \cite{BansayePardoSmadi2013}. From (\ref{eqn5.03}) and (\ref{eqn5.04}) we see that the survival probability of the CBRE-process up to time $t\ge 0$ is given by
  \begin{eqnarray}\label{eq5.3}
 \mathbf{P}(X_t(x)>0)
 \ar=\ar \mathbf{P}(Z_t(x)>0)
 =\lim_{\lambda\to \infty}\mathbf{E}\big[1-e^{-\lambda Z_t(x)}\big]
 = \mathbf{E}\Big[F_x\Big(\int_0^t e^{-\gamma \xi(s)}ds\Big)\Big],
  \end{eqnarray}
 where $F_x(z): = 1 - \exp\{-x(c\gamma z)^{-1/\gamma}\}$ satisfies Condition~\ref{ConF1} and \ref{ConF2}.
 The following theorem is an immediate consequence of Theorem~\ref{T1}. Using the notation introduced there, it gives characterizations of the thee regimes of the survival probability of the CB-process in heavy-tailed L\'evy random environment:
 \begin{theorem}\label{th5.1}
 We have the following three regimes of the survival probability:
 \begin{enumerate}
 \item[(1)] \textit{(Supercritical)} If $a<0$, then the following nonzero, finite limit exists:
     \begin{eqnarray*}
     \lim_{t\to\infty}\mathbf{P}\{X_t(x)>0\}
     =\mathbf{E}[F_x(A_\infty(\gamma\xi))\}];
     \end{eqnarray*}

 \item[(2)] \textit{(Critical)}\footnote{ \cite{BansayePardoSmadi2019} also considered this case with the general branching mechanism under an additional exponential moment condition:  $\mathbf{E}[e^{\theta^+ \xi_1}]<\infty$ for some $\theta^+>1$, which can not be satisfied by the random environment with regularly varing right tail.}
  If $a=0$ and $\rho_0:=\lim_{t\to\infty}\mathbf{P}\{\xi_t<0\}\in[0,1)$. There exists a constant $C(x)>0$ and a slowly varying function $\ell_0(x)$ at $\infty$ such that
    \begin{eqnarray*}
     \lim_{t\to\infty} t^{\rho_0}\ell_0(t)\mathbf{P}\{X_t(x)>0\}=C(x);
     \end{eqnarray*}

 \item[(3)] \textit{(Subcritical)} If $a>0$ and Condition~\ref{Con2.1} and \ref{C2.2} hold. Then
    \begin{eqnarray*}
     \lim_{t\to\infty}\frac{(at)^{\alpha}}{\ell(at)}\mathbf{P}\{X_t(x)>0\}=\int_0^\infty \mathbf{E}\big[C_{F_x}(\xi,s)\big]ds<\infty,
    \end{eqnarray*}
     where
     $
     C_{F_x}(\xi, s):= \lim_{t\to\infty}\mathbf{E}\big[{F_x}\big(A_s(\gamma\xi)+e^{-\gamma(\xi_{s}+\mathbf{J})}A_{t-s}(\gamma\hat\xi)\big)\,\big|\,\sigma(\xi),\mathbf{J}> at\big].$
 \end{enumerate}
 \end{theorem}

 \section{Appendix}
 \setcounter{equation}{0}

 \textit{Proof for Lemma~\ref{Thm109}.}
 Let $\xi'$ be an independent copy of $\xi$.
 For any $n\geq 1$,  from Proposition 4.1 in \cite{Fay2006} we have
 $n\mathbf{P}\{ \xi_{1/n} >x \} \sim  \mathbf{P}\{ \xi_{1} >x \} $ and hence for any $\delta>0$,
  \begin{eqnarray*}
   n\mathbf{P}\{ \xi_{1/n} \in( x,x+\delta ]\} = n\mathbf{P}\{ \xi_{1/n} >x \}-n\mathbf{P}\{ \xi_{1/n} >x+\delta \}\sim  \mathbf{P}\{ \xi_{1} \in (x,x+\delta] \} .
 \end{eqnarray*}
 For any $c<\epsilon/2$, from the independent increments of $\xi$
 we have $\mathbf{P}\big\{\xi_t\in[x,x+\delta)\big\}
 = \mathbf{P}\big\{\xi_{[nt]/n}+\hat\xi_{t-[nt]/n}\in[x,x+\delta)\big\}$ for any $t\geq 0$
 and the term on the left side of (\ref{eqn3.24}) can be bounded by
 \begin{eqnarray*}
 \sum_{k=1}^3 I_k(n,c,t) \ar:=\ar\sup_{x\geq (\epsilon-a)t}\int_{-ct}^{ct} \Big| \frac{\mathbf{P}\big\{\xi_{[nt]/n}+y\in[x,x+\delta)\big\}}{t\cdot \mathbf{P}\{\xi_1\in[at+x,at+x+\delta)\}}-1 \Big| \mathbf{P}\{ \hat\xi_{t-[nt]/n}\in dy \}\cr
  \ar\ar +\sup_{x\geq (\epsilon-a)t}\int_{-\infty}^{-ct} \Big| \frac{\mathbf{P}\big\{\xi_{[nt]/n}+y\in[x,x+\delta)\big\}}{t\cdot \mathbf{P}\{\xi_1\in[at+x,at+x+\delta)\}}-1 \Big| \mathbf{P}\{ \hat\xi_{t-[nt]/n}\in dy \}\cr
 \ar\ar +   \sup_{x\geq (\epsilon-a)t}\int_{ct}^\infty \Big| \frac{\mathbf{P}\big\{\xi_{[nt]/n}+y\in[x,x+\delta)\big\}}{t\cdot \mathbf{P}\{\xi_1\in[at+x,at+x+\delta)\}}-1 \Big| \mathbf{P}\{ \hat\xi_{t-[nt]/n}\in dy \}.
 \end{eqnarray*}
 It is easy to see that $I_1(n,c,t)$ can be bounded by the sum of the following two terms:
 \begin{eqnarray}\label{eqn.A01}
 \sup_{x\geq (\epsilon-a)t, \, |y|\leq ct} \left|\frac{n\mathbf{P}\{\xi_{1/n}\in[at+x-y,at+x-y+\delta)\}}{\mathbf{P}\{\xi_1\in[at+x,at+x+\delta)\}}-1\right|
\end{eqnarray}
and
 \begin{eqnarray}\label{eqn.A02}
 \ar\ar   \sup_{x\geq (\epsilon-a)t, \, |y|\leq ct}\frac{n\mathbf{P}\{\xi_{1/n}\in[at+x-y,at+x-y+\delta)\}}{\mathbf{P}\{\xi_1\in[at+x,at+x+\delta)\}}\cr
  \ar\ar
\times  \sup_{x\geq (\epsilon-a)t, \, |y|\leq ct} \Big|\frac{\mathbf{P}\{\xi_{[nt]/n} \in[x-y,x-y+\delta)\}}{nt\mathbf{P}\{\xi_{1/n}\in[at+x-y,at+x-y+\delta)\}}-1\Big|.
 \end{eqnarray}
 We first prove that the term $(\ref{eqn.A01})$ vanishes as $t\to\infty$.  Actually, it can be bounded by
  \begin{eqnarray*}
 \ar\ar \sup_{x\geq \epsilon t, \, |y|\leq ct}  \frac{\mathbf{P}\{\xi_1\geq x-y\}/(x-y)}{\mathbf{P}\{\xi_1\geq x\}/x}  \left|
 \frac{n\mathbf{P}\{\xi_{1/n}\in[x-y,x-y+\delta)\}}{\mathbf{P}\{\xi_1\geq x-y\}/(x-y)\cdot \alpha\delta}
  \frac{\mathbf{P}\{\xi_1\geq x\}/x\cdot \alpha\delta}{\mathbf{P}\{\xi_1\in[x,x+\delta)\}}-1\right| \cr
 \ar\ar +  \sup_{x\geq \epsilon t, \, |y|\leq ct} \left|\frac{\mathbf{P}\{\xi_1\geq x-y\}/(x-y)}{\mathbf{P}\{\xi_1\geq x\}/x} -1\right| .
 \end{eqnarray*}
 From Condition~\ref{Con2.1} and the fact that $x-y\geq  \frac{\epsilon}{2} \cdot t$ for any  $x\geq \epsilon t$ and $ |y|\leq ct  $, we have $\mathbf{P}\{\xi_1\geq z\}/z\sim z^{-\alpha-1 }\ell (z)  $ as $z\to \infty$ and hence there exists a constant $C>0$ independent of $n$ and $c$ such that
 \begin{eqnarray*}
 \limsup_{t\to\infty} \sup_{x\geq \epsilon t, \, |y|\leq ct} \left|\frac{\mathbf{P}\{\xi_1\geq x-y\}/(x-y)}{\mathbf{P}\{\xi_1\geq x\}/x} -1\right|  \ar\leq\ar C \Big|\Big(1-\frac{c}{\epsilon}\Big)^{-1-\alpha}-1\Big|\vee \Big|\Big(1+\frac{c}{\epsilon}\Big)^{-1-\alpha}-1\Big|,
 \end{eqnarray*}
 which vanishes as $c\to 0$.
Moreover, from Condition~\ref{C2.2} and the fact that $x-y\geq  \frac{\epsilon}{2} \cdot t$ for any  $x\geq \epsilon t$ and $ |y|\leq ct  $, we also have
  \begin{eqnarray*}
 	 \lim_{t\to\infty}\sup_{x\geq \epsilon t, \, |y|\leq ct} \frac{n\mathbf{P}\{\xi_{1/n}\in[x-y,x-y+\delta)\}}{\mathbf{P}\{\xi_1\geq x-y\}/(x-y)\cdot \alpha \delta}
 	 = \lim_{t\to\infty}\sup_{x\geq \epsilon t} \frac{\mathbf{P}\{\xi_1\in[x,x+\delta)\}}{\mathbf{P}\{\xi_1\geq x\}/x\cdot \alpha \delta}=1
 \end{eqnarray*}
 and
 \begin{eqnarray*}
 \limsup_{t\to\infty} \sup_{x\geq \epsilon t, \, |y|\leq ct}
 \left|
\frac{n\mathbf{P}\{\xi_{1/n}\in[x-y,x-y+\delta)\}}{\mathbf{P}\{\xi_1\geq x-y\}/(x-y)\cdot \alpha\delta}
\frac{\mathbf{P}\{\xi_1\geq x\}/x\cdot \alpha\delta}{\mathbf{P}\{\xi_1\in[x,x+\delta)\}}-1\right|  =0.
 \end{eqnarray*}
Taking these estimates above into (\ref{eqn.A01}), we have
 \begin{eqnarray*}
\lim_{c\to 0+}\limsup_{t\to\infty}\sup_{x\geq (\epsilon-a)t, \, |y|\leq ct} \left|\frac{n\mathbf{P}\{\xi_{1/n}\in[at+x-y,at+x-y+\delta)\}}{\mathbf{P}\{\xi_1\in[at+x,at+x+\delta)\}}-1\right| =0.
\end{eqnarray*}
We now start to  prove the term  (\ref{eqn.A02}) vanishes as $t\to\infty$. From the result above, we see that the first term in (\ref{eqn.A02}) is uniformly bounded. For the second term,   from Corollary~2.1 in \cite{DenisovDieker2008} we have as $t\to\infty$,
 \begin{eqnarray*}
 \lefteqn{\sup_{x\geq (\epsilon-a)t} \sup_{y\leq ct} \Big|\frac{\mathbf{P}\{\xi_{[nt]/n} \in[x-y,x-y+\delta)\}}{nt\mathbf{P}\{\xi_{1/n}\in[at+x-y,at+x-y+\delta)\}}-1\Big|}\ar\ar\cr
 \ar \ar\qquad \leq \sup_{x\geq (\epsilon/2-a)t} \Big|\frac{\mathbf{P}\{\xi_{[nt]/n} \in[x,x+\delta)\}}{nt\mathbf{P}\{\xi_{1/n}\in[at+x,at+x+\delta)\}}-1\Big|\to 0.
 \end{eqnarray*}
From these two results above, we have $\lim_{c\to 0+}\limsup_{t\to\infty}I_1(n,c,t)=0$.
 For $I_2(n,c,t)$, from Corollary~2.1 in \cite{DenisovDieker2008} we also have
 \begin{eqnarray*}
 I_2(n,c,t) \leq \sup_{x\geq (\epsilon+c-a)t}  \Big| \frac{\mathbf{P}\big\{\xi_{[nt]/n}\in[x,x+\delta)\big\}}{t\cdot \mathbf{P}\{\xi_1\in[at+x,at+x+\delta)\}}-1 \Big| ,
 \end{eqnarray*}
 which vanishes as $t\to \infty$.
 It is easy to see that  $I_3(n,c,t)$ can be bounded by
 \begin{eqnarray*}
  \frac{\mathbf{P}\{S_{1/n}\geq ct\}}{t\cdot \mathbf{P}\{\xi_1\in[at+x,at+x+\delta)\}}+ \mathbf{P}\{S_{1/n}\geq ct\}.
 \end{eqnarray*}
 Here the second term on the right side of the equality above goes to $0$ as $t\to\infty$. For the first term,
 from (\ref{eqn2.05}), there exists a constant $C>0$ such that
 \begin{eqnarray*}
 S_{1/n}\leq C + \sup_{s\in[0,1/n]} \Big|\sigma  B_s +\int_0^s\int_{|x|\leq 1}x\tilde{N}(ds,dx)\Big|+ \int_0^{1/n} \int_1^\infty xN(ds,dx).
 \end{eqnarray*}
 From Proposition 4.1 in \cite{Fay2006}, as $t\to\infty$ we have
 \begin{eqnarray*}
 \mathbf{P}\{S_{1}\geq ct\}\sim \mathbf{P}\Big\{  \int_0^{1/n} \int_1^\infty xN(ds,dx)\geq  ct/2 \Big\}  \ar\sim\ar \frac{\bar\nu(ct/2)}{n}
 \end{eqnarray*}
 and hence $\lim_{n\to\infty}\limsup_{t\to\infty}I_3(n,c,t)= 0$.  Here we have finished the whole proof.
 \qed

 \medskip

 \textbf{Acknowledgements.}
 The author would like to thank Professor Mladen Savov for enlightening comments.
 He also recommended the author his recent works about exponential functionals of L\'evy processes, which helped a lot to simplify the proofs in Section~4.
 The author is also grateful for the helpful comments from the two professional referees and the financial support from the Alexander-von-Humboldt-Foundation.

%
%
%
%
%
%
%
%
%
%
%

\begin{thebibliography}{}

 \bibitem[\protect\citeauthoryear{Afanasyev et al.}{2005}]{AGKV2005}
  Afanasyev, V.~I., Geiger, J., Kersting, G. and Vatutin, V. (2005) Criticality for branching processes in random environment. \textit{Ann. Probab.}, \textbf{33}(2), 645--673.


 \bibitem[\protect\citeauthoryear{Bansaye et al.}{2019}]{BansayeCaballeroMeleard2019}
 Bansaye, V., Caballero, M. E. and M\'el\'eard, S. (2019) Scaling limits of population and evolution processes in random environment. \textit{Electron. J. Probab.}, \textbf{24}, 1--38.

 \bibitem[\protect\citeauthoryear{Bansaye et al.}{2019}]{BansayePardoSmadi2019}
  Bansaye, V.,  Pardo, J.~C. and Smadi, C. (2019) Extinction rate of continuous state branching processes in critical {L}\'{e}vy environments.  \textit{arXiv:1903.06058}.
%
 \bibitem[\protect\citeauthoryear{Bansaye et al.}{2013}]{BansayePardoSmadi2013}
 Bansaye, V., Pardo~Millan, J.~C. and Smadi, C. (2013) On the extinction of continuous state branching processes with catastrophes. \textit{Electron. J. Probab.}, \textbf{106}, 1--31.
%
 \bibitem[\protect\citeauthoryear{Bansaye and Vatutin}{2017}]{BansayeVatutin2017}
  Bansaye, V. and Vatutin, V. (2017) On the survival probability for a class of subcritical branching processes in random environment. \textit{Bernoulli}, \textbf{23}(1), 58--88.

 \bibitem[\protect\citeauthoryear{Barker and Savov}{2019}]{BarkerSavov2019}
 Barker, A. and Savov, M. (2019) Bivariate {B}ernstein-gamma functions and moments of exponential functionals of subordinators. \textit{arXiv:1907.07966}.
%
 \bibitem[\protect\citeauthoryear{Bertoin}{1996}]{Bertoin1996}
  Bertoin, J. (1996) \textit{L\'evy Processes}. Cambridge University Press.
%
 \bibitem[\protect\citeauthoryear{Bertoin and Yor}{2005}]{BertoinYor2005}
  Bertoin,  J. and Yor, M. (2005) Exponential functionals of {L}\'{e}vy  processes.
   \textit{Probab. Surv.}, \textbf{2}. 191--212.

 \bibitem[\protect\citeauthoryear{Billingsley}{1999}]{Billingsley1999}
  Billingsley, P. (1999) \textit{Convergence of {P}robability {M}easures}, 2nd edn. New York: Wiley.

 \bibitem[\protect\citeauthoryear{B\"{o}inghoff et al.}{2012}]{BoinghoffHutzenthaler2012}
   B\"{o}inghoff, C. and Hutzenthaler, M. (2012) Branching diffusions in random environment. \textit{Markov Process. Relat. Fields}, \textbf{18}(2), 269--310.
%
 \bibitem[\protect\citeauthoryear{Carmona et al.}{1997}]{CarmonaPetitYor1997}
  Carmona, P, Petit, F. and Yor, M. (1997) On the distribution and asymptotic results for exponential functionals of {L}\'{e}vy processes. \textit{Exponential functionals and principal values related to Brownian motion}, 73--121.

 \bibitem[\protect\citeauthoryear{Cline}{1986}]{Cline1986}
  Cline, D.~B.~H. (1986) Convolution tails, product tails and domains of attraction. \textit{Probab. Theory Relat. Fields}, \textbf{72}(4), 529--557.

 \bibitem[\protect\citeauthoryear{Comtet et al.}{1998}]{ComtetMonthusYor1998}
 Comtet, A., Monthus, C. and Yor, M. (1998) Exponential functionals of {B}rownian motion and disordered systems.
 \textit{J. Appl. Probab.}, \textbf{35}, 255--271.

 \bibitem[\protect\citeauthoryear{Denisov et al.}{2008}]{DenisovDieker2008}
 Denisov, D.,  Dieker, A.~B. and Shneer, V. (2008) Large deviations for random walks under subexponentiality: the big-jump domain. \textit{Ann. Probab.}, \textbf{36}(5), 1946--1991.

 \bibitem[\protect\citeauthoryear{Denisov et al.}{2010}]{DenisovFossKorshunov2010}
  Denisov, D., Foss, S. and Korshunov, D. (2010) Asymptotics of randomly stopped sums in the presence of heavy tails. \textit{Bernoulli}, \textbf{16}(4), 971--994.

 \bibitem[\protect\citeauthoryear{Denisov and Shneer}{2013}]{DenisovShneer2013}
 Denisov, D. and Shneer, V. (2013) Asymptotics for the first passage times of {L}\'{e}vy processes and random walks. \textit{J. Appl. Probab.}, \textbf{50}(1), 64--84.

  \bibitem[\protect\citeauthoryear{Denisov et al.}{2014}]{DenisovVatutinWachtel2014}
 Denisov, D.,Vatutin, V. and Wachtel, V. (2014) Local probabilities for random walks with negative drift conditioned to stay nonnegative.
 \textit{Electron. J. Probab.}, \textbf{19}(88), 1-17.


 \bibitem[\protect\citeauthoryear{Doney}{1985}]{Doney1985}
 Doney, R.~A. (1985) Conditional limit theorems for asymptotically stable random walks. \textit{Probab. Theory Relat. Fields}, \textbf{70}(3), 351--360.

 \bibitem[\protect\citeauthoryear{Doney and Jones}{2012}]{DoneyJones2012}
 Doney, R.~A. and  Jones, E.~M. (2012) Conditioned random walks and {L}\'{e}vy processes. \textit{Bull. London Math. Soc.}, \textbf{44}, 139--150.
%
 \bibitem[\protect\citeauthoryear{Durrett}{1978}]{Durrett1978}
 Durrett, R. (1978)  Conditioned limit theorems for some null recurrent Markov processes. \textit{Ann. Probab.}, \textbf{6}(5), 798--828.
%
 \bibitem[\protect\citeauthoryear{Durrett}{1980}]{Durrett1980}
  Durrett, R. (1980)
  Conditioned limit theorems for random walks with negative drift.
  \textit{Probab. Theory Relat. Fields}, \textbf{52}(3), 277--287.
%
\bibitem[\protect\citeauthoryear{Fa\"{y} et al.}{2006}]{Fay2006}
  Fa\"{y}, G., Gonz\'{a}lez-Ar\'{e}valo, B., Mikosch, T. and
  Samorodnitsky, G. (2006) Modeling teletraffic arrivals by a {P}oisson cluster process. \textit{Queueing Systems}, \textbf{54}(2), 121--140.

 \bibitem[\protect\citeauthoryear{Fu and Li}{2010}]{FuLi2010}
   Fu, Z. and Li, Z. (2010) Stochastic equations of non-negative processes with jumps. \textit{Stochastic Process. Appl.}, \textit{120}(3), 306--330.

 \bibitem[\protect\citeauthoryear{Geiger et al.}{2003}]{GeigerKerstingVatutin2003}
   Geiger, J., Kersting, G. and Vatutin, V. (2003) Limit theorems for subcritical branching processes in random environment.
   \textit{Ann. Inst. H. Poincar?Probab. Statist.}, \textit{39}(4), 593--620.
%
 \bibitem[\protect\citeauthoryear{Geman and Yor}{1993}]{GemanYor1993}
   Geman, H. and Yor, M. (1993) Bessel processes, {A}sian options, and perpetuities. \textit{Math. Finance}, \textbf{3}(4), 349--375.

 \bibitem[\protect\citeauthoryear{Guivarc'h and Liu}{2001}]{GuivarchLiu2001}
   Guivarc'h, Y. and Liu, Q. (2001) Propri\'{e}t\'{e}s asymptotiques des processus de branchement en environnement al\'{e}atoire.  \textit{C. R. Acad. Sci. Paris Sér. I Math.}, \textbf{332}(4), 339--344.

 \bibitem[\protect\citeauthoryear{He et al.}{2018}]{HeLiXu2018}
   He, H., Li, Z. and Xu, W. (2018) Continuous-state branching processes in {L}\'{e}vy random environments. \textit{J. Theoret. Probab.}, \textbf{31}(4), 1952--1974.

 \bibitem[\protect\citeauthoryear{Kawazu and Tanaka}{1993}]{KawazuTanaka1993}
   Kawazu, K. and Tanaka, H. (1993) On the maximum of a diffusion process in a drifted {B}rownian environment, in \textit{S\'{e}minaire de {P}robabilit\'{e}s,  {XXVII}}, \textbf{1557}, Berlin: Springer, 78--85.

 \bibitem[\protect\citeauthoryear{Kurtz}{1978}]{Kurtz1978}
   Kurtz, T.~G. (1978) Diffusion approximations for branching processes. \textit{Branching processes, Conf. Quebec 1976, Adv. Probab. relat. Top.}.  \textit{5}, 269-292.

 \bibitem[\protect\citeauthoryear{Kyprianou}{2006}]{Kyprianou2006}
   Kyprianou, A.~E. (2006) \textit{Introductory {L}ectures on {F}luctuations of {L}\'{e}vy {P}rocesses with {A}pplications}. Berlin: Springer-Verlag.

 \bibitem[\protect\citeauthoryear{Li and Xu}{2018}]{LiXu2018}
   Li, Z. and Xu, W. (2018) Asymptotic results for exponential functionals of {L}\'{e}vy processes. \textit{Stochastic Process. Appl.}, \textbf{128}(1), 108--131.

 \bibitem[\protect\citeauthoryear{Palau and Pardo}{2017}]{PalauPardo2015}
   Palau, S. and Pardo, J.~C. (2017) Continuous state branching processes in random environment: the {B}rownian case.
   \textit{Stochastic Process. Appl.}, \textbf{127}(3), 957--994.

 \bibitem[\protect\citeauthoryear{Palau and Pardo}{2018}]{PalauPardo2018}
   Palau, S. and Pardo, J. C. (2018)  Branching processes in a {L}\'{e}vy random environment. \textit{Acta Appl. Math.}, \textbf{153}, 55--79.

 \bibitem[\protect\citeauthoryear{Palau et al.}{2016}]{PalauPardoSmadi2016}
   Palau, S., Pardo, J.~C., and Smadi, C. (2016) Asymptotic behaviour of exponential functionals of {L}\'{e}vy processes with applications to random processes in random environment. \textit{ALEA Lat. Am. J. Probab. Math. Stat.}, \textbf{13}(2), 1235--1258.

 \bibitem[\protect\citeauthoryear{Pardo et al.}{2012}]{PardoPatieSavov2012}
   Pardo, J.~C., Patie, P. and Savov, M. (2012) A {W}iener-{H}opf type factorization for the exponential functional of {L}\'evy processes. \textit{J. Lond. Math. Soc.} (2), \textbf{86}(3), 930--956.

 \bibitem[\protect\citeauthoryear{Patie and Savov}{2012}]{PatieSavov2012}
   Patie, P. and Savov, M. (2012) Extended factorizations of exponential functionals of {L}\'evy processes. \textit{Electron. J. Probab.}, \textbf{17}(38), 1--22.

 \bibitem[\protect\citeauthoryear{Patie and Savov}{2018}]{PatieSavov2018}
   Patie, P. and Savov, M. (2018) Bernstein-gamma functions and exponential functionals of {L}\'{e}vy processes. \textit{Electron. J. Probab.}, \textbf{23}(75), 1--101.

  \bibitem[\protect\citeauthoryear{Sato}{1999}]{Sato1999}
   Sato, K. (1999) \textit{L\'{e}vy processes and infinitely divisible distributions}. Cambridge University Press.

 \bibitem[\protect\citeauthoryear{Stephenson}{2018}]{Stephenson2018}
   Stephenson, R. (2018) On the exponential functional of {M}arkov additive processes, and applications to multi-type self-similar fragmentation processes and trees. \textit{ALEA, Lat. Am. J. Probab. Math. Stat.}, \textbf{15}(2), 1257--1292.
%
 \bibitem[\protect\citeauthoryear{Vatutin and Wachtel}{2009}]{VatutinWachtel2009}
   Vatutin V. and Wachtel, V. (2009)  Local probabilities for random walks conditioned to stay positive.
   \textit{Probab. Theory Relat. Fields}, \textbf{143}(1-2), 177--217.

 \bibitem[\protect\citeauthoryear{Vatutin and Zheng}{2012}]{VatutinZheng2012}
   Vatutin, V. and Zheng, X. (2012) Subcritical branching processes in a random environment without the {C}ram\'er condition. \textit{Stochastic Process. Appl.}, \textbf{122}(7), 2594--2609.

 \bibitem[\protect\citeauthoryear{V\'echambre}{2019}]{Vechambre2019}
   V\'{e}chambre, G. (2019) Exponential functionals of spectrally one-sided {L}\'{e}vy processes conditioned to stay positive.
   \textit{Ann. Inst. H. Poincar?Probab. Statist.}, \textbf{55}(2), 620--660.

 \bibitem[\protect\citeauthoryear{Yor}{1992}]{Yor1992}
   Yor, M. (1992) On some exponential functionals of brownian motion. \textit{Adv. in Appl. Probab.}, \textbf{24}(3), 509--531.
%

\end{thebibliography}
%

 \end{document}